\documentclass[11pt]{article}
\usepackage[title,titletoc]{appendix}
\usepackage{hyperref}
\hypersetup{colorlinks,allcolors=black}
\usepackage{amssymb}
\usepackage{amsmath}
\usepackage{amsthm}
\usepackage{graphicx}
\usepackage{color}
\usepackage{authblk}
\usepackage{lscape}
\usepackage{longtable}
\usepackage{soul}
\usepackage[scriptsize]{caption}
\usepackage{subcaption}
\usepackage{setspace}
\usepackage{natbib}
\usepackage{url}
\usepackage{algorithm}
\usepackage{algpseudocode}
\usepackage[left=1in,top=1in,right=1in,bottom=1in,nohead]{geometry}
\long\def\symbolfootnote[#1]#2{\begingroup%
	\def\thefootnote{\fnsymbol{footnote}}\footnote[#1]{#2}\endgroup}

\newcommand{\sumt}{\sum\limits_{t=1}^N}

\newcommand{\rhommin}{\rho_m^{\text{min}}}
\newcommand{\rhommax}{\rho_m^{\text{max}}}

\allowdisplaybreaks[4]
\begin{document}
\title{Robust fractionation in cancer radiotherapy}
\author{Ali Ajdari$^{1,2}$\footnote{Email: aajdari@mgh.harvard.edu}   , Archis Ghate$^2$}
\date{$^1$Department of Radiation Oncology, Massachusetts General Hospital \& Harvard Medical School, Boston, MA\\%
	$^2$Department of Industrial \& Systems Engineering, University of Washington, Seattle, WA\\[2ex]%
	January 2016}
\maketitle
\begin{abstract}
In cancer 
radiotherapy, the standard formulation of the optimal fractionation problem based 
on the linear-quadratic dose-response model is a
non-convex quadratically constrained quadratic program (QCQP). An optimal 
solution for this QCQP can be derived by 
solving a two-variable linear program. Feasibility of this solution, however, 
crucially depends on the so-called alpha-over-beta ratios for the organs-at-risk, 
whose true values are 
unknown. Consequently, the dosing schedule presumed optimal, in fact, may not 
even be  
feasible in practice. We address this by proposing a robust counterpart of the nominal formulation. We show that a robust solution can 
be derived by solving a small number of two-variable linear programs, each with a 
small number of constraints. We  
quantify the price of robustness, and 
compare the incidence and extent of infeasibility of the nominal and robust 
solutions via numerical experiments.
\end{abstract}
%Keywords: linear-quadratic dose-response, quadratically constrained quadratic programming, linear programming.
\section{Background and motivation}
\label{sec:intro}
The goal in external beam radiotherapy for cancer is to maximize damage to the 
tumor while limiting toxic effects of radiation on nearby organs-at-risk (OAR). 
Treatment is typically delivered over multiple treatment sessions called fractions. 
This leads to a well-known optimization problem, often referred to as the 
{\emph{fractionation problem}}. 
The goal in this problem is to find the number of fractions $N$ and a 
corresponding sequence $\vec d=(d_1,d_2,\ldots,d_N)$ of doses so as to 
maximize tumor-damage while ensuring that the OAR can safely tolerate these 
doses. The fundamental tradeoffs in this problem are as follows. Normal-cells 
often have a better damage-repair capability than tumor-cells. Temporal 
dispersion of dose across multiple fractions thus gives the OAR time to recover 
between sessions. For most tumors, a large number of fractions with a small dose 
per fraction allows the treatment planner to inflict more damage on the tumor as 
compared to administering a small number of fractions with a large dose per 
fraction. However, tumors can proliferate over the treatment course, and thus a 
shorter course might work better as it kills the tumor before any significant 
proliferation. Thus the question is whether or not and how the treatment planner 
can exploit, for patients' benefit, the differences in the way in which tumors and 
OAR respond to 
radiation. 
\subsection{Mathematical formulations of the fractionation problem}\label{sec:existingform}
The fractionation problem has been studied extensively, both 
clinically and mathematically, for over a century \cite{Rockwell:1998ve}.  Mathematical formulations of this problem routinely rely on the 
linear-quadratic (LQ) cell-survival model \cite{Hall05}. Key parameters of the LQ 
model include the so-called $\alpha/\beta$ ratios for the
OAR. Research that uses this LQ model has evolved from 
single-OAR formulations, to two-OAR formulations, and, more recently, to models 
with multiple OAR. All of these formulations belong to the class of 
non-convex quadratically constrained quadratic programs (QCQPs) --- problems known to be computationally difficult in general. A closed-form optimal 
solution is available for the single OAR 
case (see, for example, 
\cite{bortfeld13,fowler2008,Fowler1995,Jones1995,Mizuta2012,unkelbach13} and 
references therein). One paper provided an optimal dosing 
scheme using Karush-Kuhn-Tucker conditions for the two-OAR 
case for a fixed $N$ \cite{bertuzzi2013}. A simulated annealing heuristic was 
applied to a two-OAR formulation in \cite{Yang2005fractionation}. 

The most recent multiple-OAR formulation of this problem (see 
\cite{Saberian2015mathmedbio,Saberian2015}) is given by
\begin{align}
\label{optfracobj}
{\text{(FRAC)}}\ \max\limits_{\vec d,\ N}\ \ &\alpha_0\sumt d_t+\beta_0 \sumt d^2_t-\tau(N)\\
\label{optfraccon} \sumt d_t+\rho_m \sumt d_t^2 &\leq \text{BED}_m,\ m\in\mathcal M,\\
\label{optfracconnneg}\vec d&\geq 0,\\
\label{optfracNcon} 1&\leq N \leq N_{\text{max}},\ \text{integer}.
\end{align}
In this problem, $\alpha_0, \beta_0$ are the tumor's 
dose-response parameters as per the LQ model. The term 
$\tau(N)$ in the objective function accounts for tumor proliferation and is given by
\begin{equation}
\label{eqn:tumorproliferation}
\tau(N)=\frac{\left[(N-1)-T_{\text{lag}}\right]^+{\ln}\ 2}{T_{\text{double}}},
\end{equation}
where $\left[(N-1)-T_{\text{lag}}\right]^+$ is defined as $\max\left\{0,(N-1)-
T_{\text{lag}}\right\}$. Here, $T_{\text{lag}}$ is the time-lag (in days) after 
which tumor proliferation 
starts after 
treatment initiation; and $T_{\text{double}}$ (in days) denotes the doubling time 
for the 
tumor. This proliferation term assumes that a single fraction is administered every 
day; it can be generalized to accommodate other fractionation schemes as 
described in \cite{Saberian2015mathmedbio}. The objective function equals the biological effect (BE) of $\vec d$ on the tumor, which is to be maximized. In constraints (\ref{optfraccon}), 
$\mathcal M=\{1,2,\ldots,n\}$ is the set of $n\geq 1$ OAR. The 
parameter $\rho_m=\beta_m/\alpha_m$ is the aforementioned (inverse) ratio of dose-response parameters for OAR $m\in \mathcal M$. The left hand side of each constraint equals the biologically effective dose (BED) administered to the corresponding OAR. The term on the right hand side is given by $\text{BED}_m=D_m+\rho_m D^2_m/N_m$. It equals the BED of a conventional treatment schedule that administers a total dose of $D_m$ in $N_m$ 
equal-dosage fractions and that OAR $m$ is known to tolerate. Thus, each of 
these constraints ensures that, for each OAR, the BED of $\vec d$ is no more than what is safe for that OAR. In constraint 
(\ref{optfracNcon}), $N_{\text{max}}$ is the maximum number of fractions 
that is logistically feasible in the treatment protocol. In the sequel, we will often refer to (FRAC) as the nominal problem.
\subsection{Optimal solution of the nominal fractionation problem}\label{sec:existingoptsol}
An optimal 
solution for this multiple-OAR case was provided in \cite{Saberian2015mathmedbio}; this
solution works either when $\alpha_0/\beta_0\leq \min\limits_{m\in \mathcal M}\  (\alpha_m/\beta_m)$ or 
when $\alpha_0/\beta_0\geq \max\limits_{m\in \mathcal M}\ (\alpha_m/\beta_m)$. 
The first provably optimal solution that works irrespective of the ordering of these 
ratios for the multiple-OAR case was recently derived in \cite{Saberian2015} 
based on the 
doctoral dissertation of Saberian \cite{fatemehthesis}. This solution 
was obtained by equivalently reformulating 
(FRAC) for each fixed $N$ as a {\emph{two-variable}} linear program (LP) with $n
+2$ linear constraints and non-negativity constraints on the two variables. 
The two variables in this LP are $x=\sumt d_t$ and $y=
\sumt d^2_t$ and the LP is given by
\begin{align}
\text{(2VARLP)}\ \max\limits_{x,y}\ \alpha_0x&+\beta_0 y-\tau(N)\\
x+\rho_m y &\leq \text{BED}_m,\ m\in\mathcal M,\\
y &\leq \gamma^* x,\\
c^*x&\leq y,\\
x&\geq 0,\\
y&\geq 0,
\end{align}
where $\gamma^*= \min\limits_{m\in\mathcal M}  b_m(1)$ and $c^*= \min\limits_{m\in\mathcal M}  b_m(N)$ with 
$b_m(N)=\frac{-1+\sqrt{1+4\rho_m\text{BED}_m/N}}{2\rho_m}$ for $m\in \mathcal M$ and for all $N\geq 1$. Specifically,
for each fixed $N$, if $x^*$, $y^*$ is an optimal solution of this LP, then the dosing 
schedule $(q,\underbrace{p,p,\ldots,p}_{N-1\ \text{times}})$, where
\begin{align}
\label{eqn:p}
p&=\frac{x^*}{N}\Bigg[1-\sqrt{1-\Bigg(1-\frac{y^*}{(x^*)^2}\Bigg)\Bigg(\frac{N}{N-1}\Bigg)}\Bigg],\\
\label{eqn:q}
q&=x^*-(N-1)p,
\end{align}
is optimal. Moreover, it can be shown that there are only three possibilities for $x^*$ and $y^*$. The first is where $\sqrt{y^*}=x^*$ and then $p=0$ (this is called a single-dosage solution); the second is where $\sqrt{Ny^*}=x^*$ and then $p=q$ (this is called an equal-dosage solution); and the third is where $\sqrt{y^*}<x^*<\sqrt{Ny^*}$ and then $0\neq p\neq q\neq 0$ (this is called an unequal-dosage solution) (see \cite{fatemehthesis,Saberian2015} for details). An optimal number 
of fractions can then be found by substituting a dosing schedule so obtained 
into the objective function in (FRAC) for each $N\in \{1,2,\ldots,N_{\text{max}}\}$ 
and picking the one that yields the largest tumor BE. Consequently, (FRAC) is 
solved by solving exactly $N_{\text{max}}$ two-variable LPs.
\subsection{Limitations of existing formulations and our contributions}\label{sec:existinglimit}
One drawback of all aforementioned formulations of the fractionation problem 
based on the LQ model is that the values of $\rho_m$ are not 
known. Thus, a dosing schedule derived using estimated or ``nominal" values of 
these parameters may not even be feasible in practice. 

In a recent unpublished manuscript \cite{badri2015v1}, Badri et al., independently of an earlier (May 2015) unpublished version of our present work, attempted to remedy this by 
studying a robust 
formulation of the above fractionation problem. In their formulation, the treatment 
planner derives a robust solution by assuming that the $\rho_m$ 
values vary within a known non-negative interval. However, the crucial dependence 
of the right hand side $\text{BED}_m$ on 
$\rho_m$ in constraints (\ref{optfraccon}) was ignored in that manuscript. This 
meant that an optimal solution to their robust formulation was obtained by 
replacing $\rho_m$ on the left hand side in (\ref{optfraccon}) by its largest 
possible value. This implied that the robust solution is derived simply by solving 
the two-variable LP in \cite{fatemehthesis,Saberian2015}. Unfortunately, since the 
right hand side in constraints (\ref{optfraccon}) in fact explicitly depends on $
\rho_m$, such a simplified solution might not be robust in practice. Badri et al. 
rectified this limitation in an updated unpublished variation \cite{badri2015v2} of 
their original manuscript, again independently of the earlier (May 2015)
unpublished version of our present work that they cited. 

The main focus of the original and the updated versions by Badri et al. was on a 
chance constrained formulation of the problem, which required the treatment planner to know the probability distribution of alpha-over-beta ratios, and which 
called for a computationally more demanding solution approach than what is 
needed for the robust formulation. On the plus side, a potential benefit of the 
resulting chance constrained solution is that it might be less conservative than the 
robust solution (although this 
is perhaps impossible to verify rigorously). Given their 
alternative focus, Badri et al. gave a somewhat cursory treatment to the robust 
approach in both their manuscripts, did not present an infeasibility analysis of the 
resulting robust solutions, did not quantify the price of robustness, and 
only included minimal sensitivity results.

Here we study essentially the same robust problem as in the updated 
version of Badri et al. We do, however, provide mathematical and clinical insights 
missing in their work. Firstly, we present our solution approach in much more 
detail. We 
show that, for each fixed $N$, an optimal solution to the non-convex robust 
problem can be recovered by solving $n+1$ two-variable LPs. 
Consequently, the robust fractionation problem is solved by solving  
$(n+1)N_{\text{max}}$ two-variable LPs; each of these LPs includes $n+2$ linear 
constraints and non-negativity constraints on the two variables. We perform  
sensitivity analyses 
with respect to the values of $T_{\text{lag}}$ and $T_{\text{double}}$ currently 
available in the clinical literature to numerically quantify the price of robustness. 
We also provide qualitative and quantitative comparisons between the nominal 
and robust fractionation schedules. Finally, we present an extensive analysis of 
the infeasibility suffered by the nominal and robust solutions in a broad range of 
scenarios.

This paper is organized as follows. Our robust formulation is described in the next 
section. The solution approach is detailed in Section \ref{sec:solve}. Numerical 
results are presented in Section \ref{sec:results}. We conclude with a summary of 
our contributions, an outline of some variations and limitations of our model, and 
opportunities for future work.
\section{A robust formulation}\label{sec:robust}
We refer the reader to \cite{bental} for a textbook and to 
\cite{robustreview} for a survey on robust optimization. We employ a standard 
interval uncertainty model from these existing works to construct a robust counterpart of the nominal problem (FRAC). Specifically, we 
use $\tilde\rho_m$ to denote the ``true" unknown value of $\rho_m$, for 
$m=1,\ldots,n$. We assume that this unknown value belongs to a 
known interval of values $[\rhommin,\rhommax]$; here $0<\rhommin\leq 
\rhommax<\infty$.
We wish to find an $N, \vec d$ pair that is feasible to BED constraints 
(\ref{optfraccon}) for all $m\in \mathcal M$ no matter what true values 
$\tilde\rho_m$ are realized (as long as they belong to the aforementioned intervals). The resulting robust counterpart of 
(FRAC) is given by
\begin{align}
\label{robj}
\max_{\vec{d},\ N}\ \ & \alpha_0 \sum \limits_{t=1}^N d_t + \beta_0\sum \limits_{t=1}^N d_t^2-\tau(N)\\
\label{rcon1}
\sumt d_t+\tilde{\rho}_m \Bigg(\sumt d_t^2-\frac{D_m^2}{N_m}\Bigg) &\leq D_m,\; \; m \in \mathcal{M},\;\forall \tilde{\rho}_m \in [\rho_m^{\text{min}},\rho_m^{\text{max}}],\\
\label{rcon2}\vec d&\geq 0,\\
\label{rcon3} 1&\leq N \leq N_{\text{max}},\ \text{integer}.
\end{align}
Note here that, for simplicity of exposition, our formulation does not consider 
uncertainty in the values of $\alpha_0$ and $\beta_0$ for the tumor. It is standard 
in robust optimization to not include uncertainty in the objective function 
coefficients. Uncertainty in these tumor 
parameters can, however, be easily incorporated by maximizing the worst-case 
value of the 
objective function (we accomplish this in our numerical results in 
Section \ref{sec:tumoruncertain}). 

%To write a compact robust formulation, it is convenient to normalize the intervals $
%[\rho_m^{\text{min}},\rho_m^{\text{max}}]$, for $m\in \mathcal M$. That is, we 
%express $\tilde \rho_m$ as $\rho_m^{\text{mean}}+\eta_m \rho_m^{\text{range}}$,  
%where $\eta_m \in [-1,1]$ is a real number; 
%$\rho_m^{\text{mean}}=(\rho_m^{\text{max}}+\rho_m^{\text{min}})/2$ and $
%\rho_m^{\text{range}}=(\rho_m^{\text{max}}-\rho_m^{\text{min}})/2$.
By introducing $\rho_m^{\text{mean}}=(\rho_m^{\text{max}}+\rho_m^{\text{min}})/2$ and $\rho_m^{\text{range}}=(\rho_m^{\text{max}}-\rho_m^{\text{min}})/2$, and after some simple algebra, we can see that for each OAR $m \in \mathcal M$, constraint \ref{rcon1} is equivalent to the following constraint: 
\begin{gather}
\sumt d_t+ \rho_m^{\text{mean}} \sumt d_t^2 + 
\rho_m^{\text{range}} \left|\sumt d_t^2-
\frac{D_m^2}{N_m} \right| \leq D_m+\rho_m^{\text{mean}}\frac{D_m^2}
{N_m}.
\end{gather}
%\begin{align*}
%&\sumt d_t +\tilde \rho_m \Bigg(\sumt d_t^2-\frac{D_m^2}{N_m}\Bigg) \leq D_m, \ 
%\forall \tilde \rho_m\in[\rho_m^{\text{min}},\rho_m^{\text{max}}]\\
%&\Leftrightarrow \sumt d_t+(\rho_m^{\text{mean}}+\eta_m \rho_m^{\text{range}}) 
%\Bigg(\sumt d_t^2-\frac{D_m^2}{N_m}\Bigg) \leq D_m,\ \forall \eta_m\in[-1,1]\\
%&\Leftrightarrow\sumt d_t+\rho_m^{\text{mean}} \Bigg(\sumt d_t^2-\frac{D_m^2}
%{N_m}\Bigg) +\eta_m 
%\rho_m^{\text{range}} \Bigg(\sumt d_t^2-\frac{D_m^2}{N_m}\Bigg) \leq D_m,\ 
%\forall \eta_m\in[-1,1]\\
%&\Leftrightarrow \sumt d_t+\rho_m^{\text{mean}} \Bigg(\sumt d_t^2-\frac{D_m^2}
%{N_m}\Bigg) +\rho_m^{\text{range}} \left| \sumt d_t^2-\frac{D_m^2}{N_m} \right| 
%\leq D_m\\
%&\Leftrightarrow\sumt d_t+ \rho_m^{\text{mean}} \sumt d_t^2 + 
%\rho_m^{\text{range}} \left|\sumt d_t^2-
%\frac{D_m^2}{N_m} \right| \leq D_m+\rho_m^{\text{mean}}\frac{D_m^2}
%{N_m}.
%\end{align*}
Thus, by defining the shorthand notation $\text{RC}_m=D_m+
\rho_m^{\text{mean}}\frac{D_m^2}{N_m}$, and putting the above pieces together, 
we can rewrite the robust counterpart (\ref{robj})-(\ref{rcon3}) as
\begin{align}
\label{rfracobj}
\text{(RFRAC)}\ f^*=\max_{\vec{d},\ N}\ \ & \alpha_0\sumt d_t+\beta_0 \sumt  d_t^2-\tau(N)
\\
\label{rfraccon1}
\sumt d_t+  \rho_m^{\text{mean}} \sumt d_t^2 + \rho_m^{\text{range}} \left| \sumt  d_t^2-\frac{D_m^2}{N_m} \right| &\leq \text{RC}_m, \ m \in \mathcal{M},
\\
\label{rfraccon2}
\vec{d} &\geq 0,\\
\label{rfraccon3}
1 &\leq  N\leq N_{\text{max}},\ \text{integer.}
\end{align}

As in the nominal problem, in order to solve this robust problem, we first solve the problems obtained by fixing $N$ 
at $1,2, \ldots, N_{\text{max}}$. For each fixed $N$, let $\vec d^*(N)=(d^*_1(N),
\ldots,d^*_N(N))$ denote the corresponding optimal dosing sequence. We then 
compare the objective values of these $N$ dosing sequences and pick 
the best. Thus, the problem we need to solve for each fixed 
$N\in\{1,2,\ldots,N_{\text{max}}\}$ is given by
\begin{align}
\label{rfracNobj}
\text{(RFRAC(N))} \ f^*(N)=\max_{\vec{d}}\  \ & \alpha_0\sumt d_t+\beta_0 
\sumt  d_t^2-\tau(N)
\\
\label{rfracNcon1}
\sumt d_t+  \rho_m^{\text{mean}} \sumt d_t^2 + \rho_m^{\text{range}} \left| \sumt  d_t^2-\frac{D_m^2}{N_m} \right| &\leq \text{RC}_m, \ m \in \mathcal{M},
\\
\label{rfracNcon2}
\vec{d} &\geq 0.
\end{align}
Note that when $\rhommin=\rhommax=\rho_m$, for 
$m=1,2,\ldots,n$, that is, when there is no uncertainty in these dose-response parameters, (RFRAC(N)) reduces to the nominal QCQP (FRAC) with $N$ 
fixed as presented in Section \ref{sec:intro}, and which was solved recently as a two-variable LP in \cite{fatemehthesis,Saberian2015}. Note, however, that the objective 
function as well as the constraints in (RFRAC(N)) are non-convex, and the 
problem is at least as hard as the nominal QCQP. 
The objective function in the 
nominal QCQP is identical in form to what we have in (RFRAC(N)), but the 
convex, quadratic constraints in the nominal QCQP do not include the absolute 
value term that appears in the corresponding constraints in (RFRAC(N)). Specifically, it is this absolute value term that makes the robust 
counterpart harder to solve as compared to the nominal problem. To overcome 
this challenge, we decompose the feasible region of (RFRAC(N)) into $n+1$
subregions in a way such that the problem over each subregion can be solved via 
a two-variable LP. The details of this procedure are discussed in the next section.
\section{Optimal solution of the robust formulation}\label{sec:solve}
To handle the absolute value term on the left hand side in constraints 
(\ref{rfracNcon1}), we decompose the non-negative orthant $\{\vec d\in\Re^N|\vec d\geq 0\}$ as 
follows. For each OAR $m\in\mathcal M$, consider two possibilities: the first is 
where $\sumt d^2_t\geq D^2_m/N_m$ and the second is $\sumt d^2_t <D^2_m/
N_m$. Suppose, in the rest of this section, without loss of generality that $D^2_1/
N_1\leq D^2_2/N_2\leq \ldots \leq D^2_n/N_n$. Then, if there is a $\vec d\geq 0$ 
and an OAR $m\in\mathcal M$ such that $\sumt d^2_t\geq D^2_m/N_m$, then for 
this $\vec d$, we have that $\sumt d^2_t\geq D^2_{m'}/N_{m'}$ for all ${m'}<m$. Similarly, if there is a $\vec d\geq 0$ 
and an OAR $m\in\mathcal M$ such that $\sumt d^2_t< D^2_m/N_m$, then for 
this $\vec d$, we have that $\sumt d^2_t < D^2_{m'}/N_{m'}$ for all ${m'}>m$. This 
means that the non-negative orthant $\{\vec d\in \Re^N| \vec d\geq 0\}$ is 
partitioned into $n+1$ subregions indexed by $k=0,1,2,\ldots,n
$. In the $k$th region, $\sumt d^2_t\geq D^2_m/N_m$ for the 
{\bf{\emph{first}}} $k$ OAR and $\sumt d^2_t<D^2_m/N_m$ for the 
{\bf{\emph{last}}} $n-k$ OAR. Let $\text{RC}^+=D_m+\rho_m^{\text{max}}
\frac{D^2_m}{N_m}$ and $\text{RC}^-=D_m+\rho_m^{\text{min}}\frac{D^2_m}
{N_m}$. Then, simple algebra reveals that for all $\vec d$ in the $k$th subregion, 
constraint (\ref{rfracNcon1}) reduces to $\sumt d_t+\rho_{m}^{\text{max}} \sumt  
d_t^2 \leq \text{RC}^+_{m}$ for OAR $m=1,2,\ldots,k$ when $k\neq 0$; and it reduces to 
$\sumt d_t+\rho_{m}^{\text{min}} \sumt  
d_t^2 \leq \text{RC}^-_{m}$ for OAR $k+1, k+2, \ldots, n$ when $k\neq n$.

As a result of the above discussion, (RFRAC(N)) is solved by solving $n+1$ 
subproblems and 
then picking a dosing schedule with the largest tumor BE from the resulting $n
+1$ solutions. The $k$th subproblem is given by
\begin{align}
\label{subkobj}
\text{(kSub(N))}\  \max_{\vec{d}}\ \ & \alpha_0\sumt d_t+\beta_0\sumt  d_t^2-\tau(N)
\\
\label{subkcon1plus}
\sumt d_t+\rho_{m}^{\text{max}} \sumt  
d_t^2 &\leq \text{RC}^+_{m},\ m=1,2,\ldots,k,\ k\neq 0,\\ 
\label{subkcon1minus}
\sumt d_t+\rho_{m}^{\text{min}} \sumt  
d_t^2 &\leq \text{RC}^-_{m},\ m=k+1,k+2,\ldots,n,\ k\neq n,\\ 
\label{subkcon2plus}
\sumt d^2_t & \geq \frac{D^2_m}{N_m},\ m=1,2,\ldots,k,\ k\neq 0,\\ 
\label{subkcon2minus}
\sumt d^2_t & < \frac{D^2_m}{N_m},\ m=k+1,k+2,\ldots,n,\ k\neq n,\\
\label{subkcon3}
\vec{d} &\geq 0.
\end{align}
In addition, owing to the fact that $D^2_1/
N_1\leq D^2_2/N_2\leq \ldots \leq D^2_n/N_n$, the group of $n$ constraints in 
(\ref{subkcon2plus})-(\ref{subkcon2minus}) reduces to at most two constraints: $
\sumt d^2_t \geq \frac{D^2_k}{N_k}$ when $k\neq 0$ and $\sumt d^2_t < 
\frac{D^2_{k+1}}{N_{k+1}}$ when $k\neq n$. After replacing this second strict 
inequality with a non-strict inequality\footnote{This can be rigorously justified by 
proving that if there is a feasible dosing schedule that satisfies $\sumt d^2_t=
\frac{D^2_{k+1}}{N_{k+1}}$ in the $k$th subproblem with a non-strict inequality, 
then this dosing schedule is feasible to the $k+1$st subproblem with a strict 
inequality; consequently, using non-strict inequalities does not alter optimality in 
our overall group of $n+1$ subproblems with strict inequalities. We omit the details of this proof for brevity.}, this simplifies the $k$th subproblem to
\begin{align}
\label{simplesubkobj}
\text{(kSub(N))}\  f^*(N;k)=\max_{\vec{d}}\ \ & \alpha_0\sumt d_t+\beta_0 \sumt  d_t^2-\tau(N)
\\
\label{simplesubkcon1plus}
\sumt d_t+\rho_{m}^{\text{max}} \sumt  
d_t^2 &\leq \text{RC}^+_{m},\ m=1,2,\ldots,k,\ k\neq 0,\\ 
\label{simplesubkcon1minus}
\sumt d_t+\rho_{m}^{\text{min}} \sumt  
d_t^2 &\leq \text{RC}^-_{m},\ m=k+1,k+2,\ldots,n,\ k\neq n,\\ 
\label{simplesubkcon2plus}
\sumt d^2_t & \geq \frac{D^2_k}{N_k},\ k\neq 0,\\ 
\label{simplesubkcon2minus}
\sumt d^2_t & \leq \frac{D^2_{k+1}}{N_{k+1}},\ k\neq n,\\
\label{simplesubkcon3}
\vec{d} &\geq 0.
\end{align}

\begin{figure}[!h]
	\centering
	\includegraphics[width=7in]{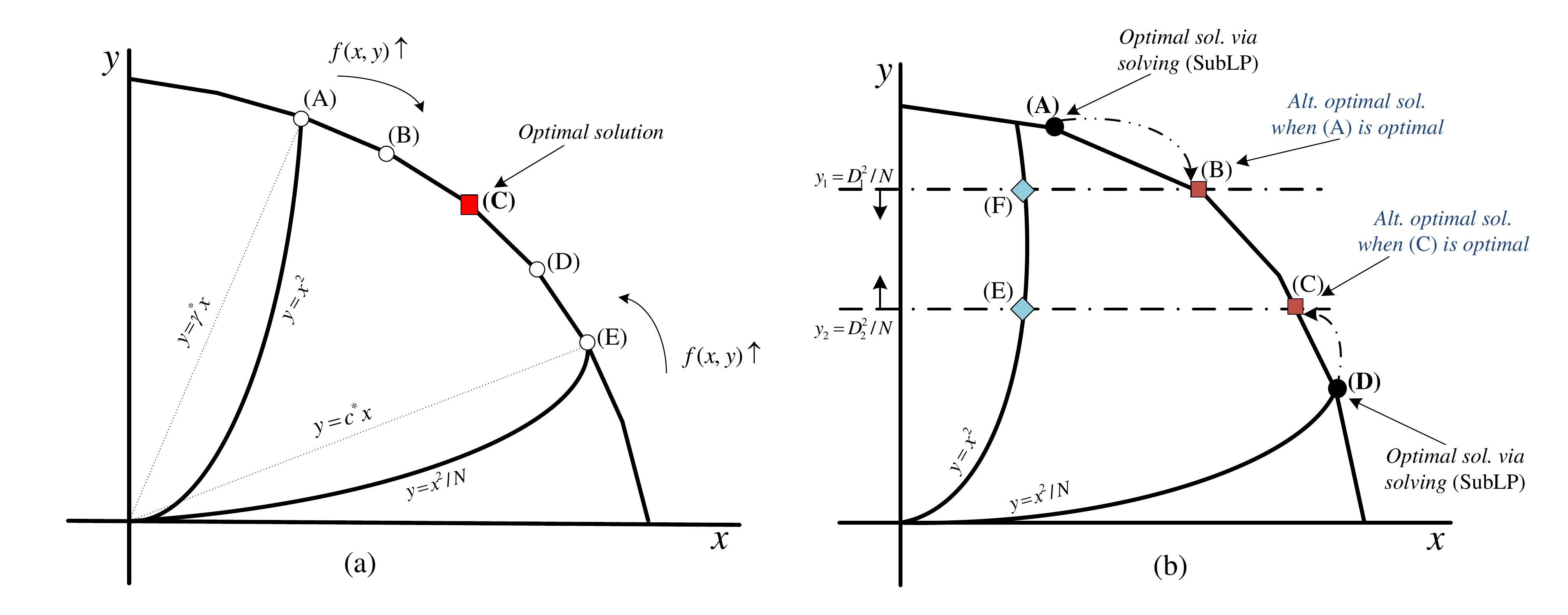}
	\caption{(a) A schematic illustration of the relationship between adjacent points in the feasible region of the subproblem $k$.  When (C) is optimal, due to the linearity of the constraint and the objective function, it can be easily shown (via geometric or analytical proof) that the denoted points in the left-hand and right-hand side of the optimal point (C) are ranked as follows in terms of their objective value: (A)$\preceq$(B)$\preceq$(C) and similarly, (E)$\preceq$(D)$\preceq$(C). (b) A schematic illustration of the method for deriving the alternative optimal solution when optimal solution obtained by solving the relaxed subproblem (SubLP $k$) does not belong to the feasible region of subproblem $k$. By following the same logic as in (a), we can argue that when (A) or (D) are the primary optimal solutions and are cut off by $y_1$ or $y_2$, (B) or (C) are the next-best optimal solutions, respectively. Also, because of the slope of the objective function, (B) and (C) are superior to (F) and (E), respectively.}
	\label{fig:multiconstfeasreg}
\end{figure}

The objective function and the constraints 
(\ref{simplesubkcon1plus})-(\ref{simplesubkcon1minus}) in this subproblem are 
identical in form to that in the nominal problem (FRAC) with $N$ fixed. Thus, the 
only difference between this subproblem and the nominal problem is the 
appearance of the additional constraints 
(\ref{simplesubkcon2plus})-(\ref{simplesubkcon2minus}). In order to solve this problem, we first relax these two constraints and use the variable transformation $x=\sumt d_t$ and $y=\sumt d^2_t$ as in \cite{fatemehthesis,Saberian2015}, to convert the relaxed 
subproblem into an equivalent two-variable LP. To write this LP compactly, we first introduce additional notation. Let
\begin{equation}
\text{RC}^k_m=
\begin{cases}
\text{RC}^+_m,\ \text{for}\ m=1,2,\ldots,k,\ k\neq 0,\\
\text{RC}^-_m,\ \text{for}\ m=k+1,k+2,\ldots,n,\ k\neq n;
\end{cases}
\end{equation}
and similarly, 
\begin{equation}
\rho^k_m=
\begin{cases}
\rho^{\text{max}}_m,\ \text{for}\ m=1,2,\ldots,k,\ k\neq 0,\\
\rho^{\text{min}}_m,\ \text{for}\ m=k+1,k+2,\ldots,n,\ k\neq n.
\end{cases}
\end{equation}
Moreover, let
\begin{align}
c_k&=\min\limits_{m\in\mathcal M}\ \frac{-1+\sqrt{1+4\rho^k_m\text{RC}^k_m/N}}{2\rho^k_m},\ \text{and}\\
\gamma_k&=\min\limits_{m\in\mathcal M}\ \frac{-1+\sqrt{1+4\rho^k_m\text{RC}^k_m}}{2\rho^k_m}.
\end{align}
Then, the two-variable LP can be written as
\begin{align}
\label{simplesubkobj}
\text{(2VARLPkSub(N))}\ \ \max_{x, y}\ \ & \alpha_0 x+\beta_0 y-\tau(N)
\\
\label{LPsimplesubkcon1}
x+\rho^k_{m} y &\leq \text{RC}^k_{m},\ m\in\mathcal M,\\ 
\label{LPsimplesubkcon2plus}
%y & \geq \frac{D^2_k}{N_k},\ k\neq 0,\\ 
%\label{LPsimplesubkcon2minus}
%y & \leq \frac{D^2_{k+1}}{N_{k+1}},\ k\neq n,\\
%\label{LPsimplesubkcon3}
y&\leq \gamma_k x,\\
\label{LPsimplesubkcon4}
c_k x &\leq y,\\
\label{LPsimplesubkcon5}
x &\geq 0,\\
\label{LPsimplesubkcon6}
y &\geq 0.
\end{align}
Let $x^*, y^*$ be an optimal solution to this two-variable LP. 
If $\frac{D^2_k}{N_k}\leq y^* \leq \frac{D^2_{k+1}}{N_{k+1}}$ as required by 
constraints (\ref{simplesubkcon2plus})-(\ref{simplesubkcon2minus}), we are done. If not, then an optimal solution can be recovered as explained in Figure \ref{fig:multiconstfeasreg}. Then, finally, a corresponding dosing schedule 
$\vec d^*(N)=(d^*_1(N),d^*_2(N),\ldots,d^*_N(N))=(q,\underbrace{p,p,\ldots,p}
_{N-1\ \text{times}})$ that is optimal to problem (kSub(N)) is recovered by 
formulas (\ref{eqn:p})-(\ref{eqn:q}).

\section{Numerical experiments}\label{sec:results}
\subsection{Qualitative properties of robust solutions}\label{sec:qualrobust}
An unavoidable downside to using robust optimization in general is that it 
sacrifices the value of the objective function in favor of a robust solution. Thus, 
it is important to quantify how much we are losing in terms of the objective value 
by solving the robust problem instead of solving the nominal problem. This is often called the price of robustness.

In this section, we numerically quantify the price of robustness via computer 
simulations for head-and-neck cancer. In these simulations, as in \cite{Saberian2015mathmedbio}, we considered four 
OAR ($n=4$), namely, spinal cord, brainstem, left and right parotids. For the 
tumor, we fixed $\alpha_0=0.35\ \text{Gy}^{-1}$ and $\beta_0=0.035\ \text{Gy}^{-2}$ 
as is standard in the clinical literature 
\cite{fowler1990,fowler01,fowler2007,fowler2008}. Again, based on the clinical 
literature \cite{fowler1990,fowler2007,Yang2005fractionation}, the nominal $
\alpha/\beta$ ratios for spinal cord, brainstem, left and 
right parotids were fixed at $3,4,5,6\ \text{Gy}$; that is, $\rho_1=1/3,\rho_2=1/4,\rho_3=1/5,\rho_4=1/6$. Here, we used different values of nominal ratios for different OAR
to fully explore the various possibilities that could arise in a robust formulation with 
multiple OAR. The tolerance doses for these OAR were fixed at  
45, 50, 26, and 28 Gy, respectively, and the conventional number of 
fractions $N_m$ was fixed at 35 days for all OAR similar to the standard QUANTEC treatment 
protocol \cite{QUANTEC}. $N_{\text{max}}$ was set to 
100 days. The uncertainty intervals were parameterized as $\tilde\rho_m\in [(1-\delta) 
\rho_m, (1+\delta) \rho_m]$, where $\delta \in [0,1]$. This allowed us to easily 
quantify the price of robustness as a function of the uncertainty level $\delta$.
We varied $\delta$ 
from 0 (to represent the nominal case) to 1 (to denote the most uncertain case 
with 
100\% uncertainty) in increments of 0.1. All experiments were carried out in 
MATLAB 
on a laptop with 2.20 GHz Intel Core2 Duo CPU and 2 GB of memory, running a 
Microsoft Windows 8.1 operating system. Tables \ref{table:PriceofRobustness} and 
\ref{table:comparisonN&d} summarize the results of our experiments for different 
values of $T_{\text{lag}}$, $T_{\text{double}}$, and $\delta$. In these tables, 
$T_{\text{lag}}$ values were set to $7, 14, 21, 28, 35$ days based on 
\cite{fowler2007} and $T_{\text{double}}$ values were set to $2,8,10,20,40,50,80, 
100$ days based on 
\cite{fowler1990,fowler2007,qi2012,Yang2005fractionation}. We are aware that 
the value of 35 days for $T_{\text{lag}}$ is perhaps too high; similarly, the values 
of 80 and 100 days for $T_{\text{double}}$ are also perhaps too high for 
head-and-neck cancer. These somewhat extreme values were included in our 
simulations to fully explore possible trends in various results of interest.

Table \ref{table:PriceofRobustness} shows, as expected, that the price of 
robustness increases with increasing $\delta$ for each $T_{\text{lag}}$, 
$T_{\text{double}}$ combination. Overall, the price of robustness seems to be 
quite small in most experiments with an average of 1.27\% over all 400 
experiments. The first, second, and third 
quartiles were 0.12\%, 0.47\%, and 1.44\%, respectively.

\begin{figure}
	\centering
	\begin{subfigure}[b]{0.5\textwidth}
		\includegraphics[width=\textwidth]{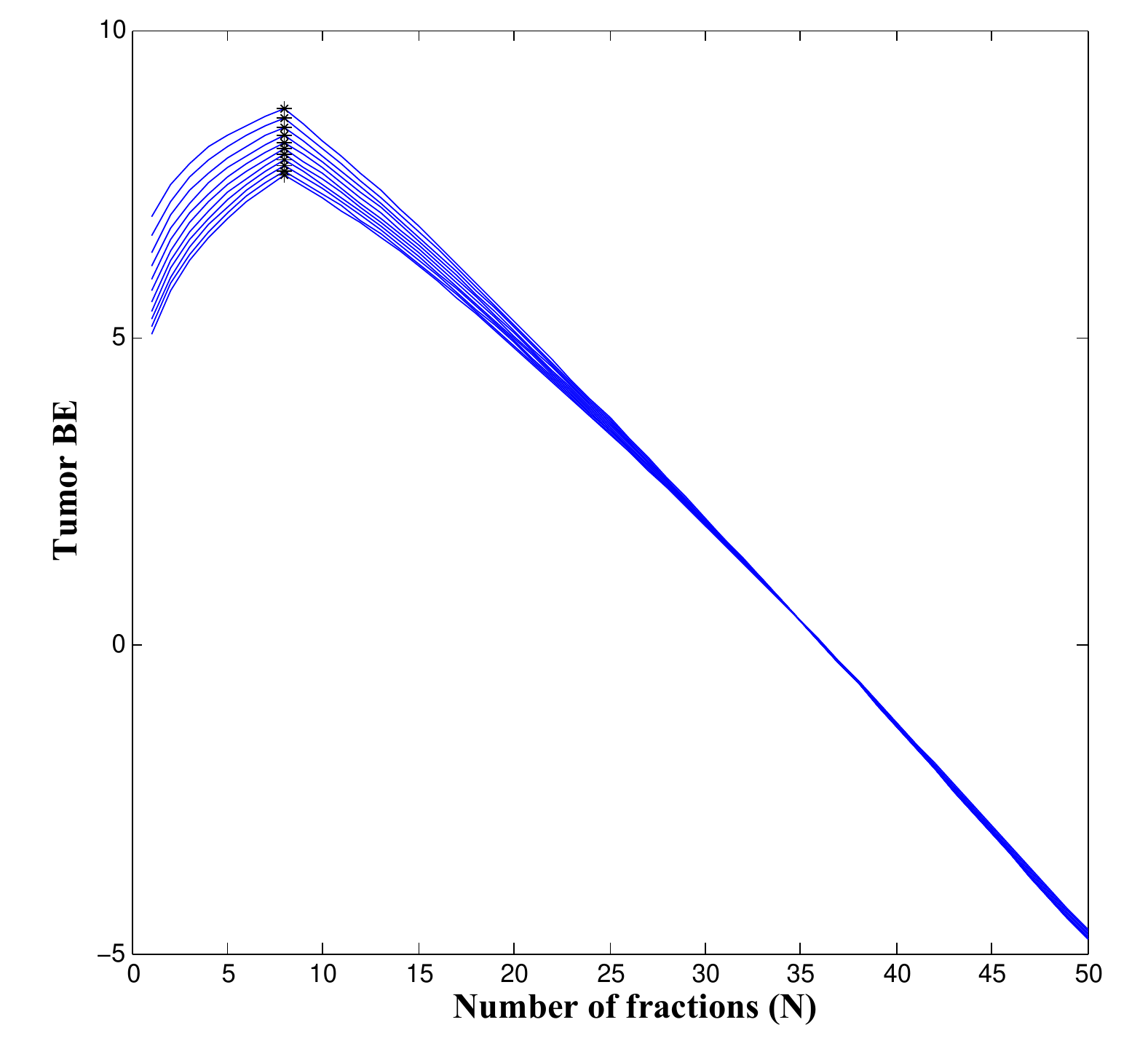}
		\caption{}
		\label{fig:7-2}
	\end{subfigure}%
	\begin{subfigure}[b]{0.5\textwidth}
		\includegraphics[width=\textwidth]{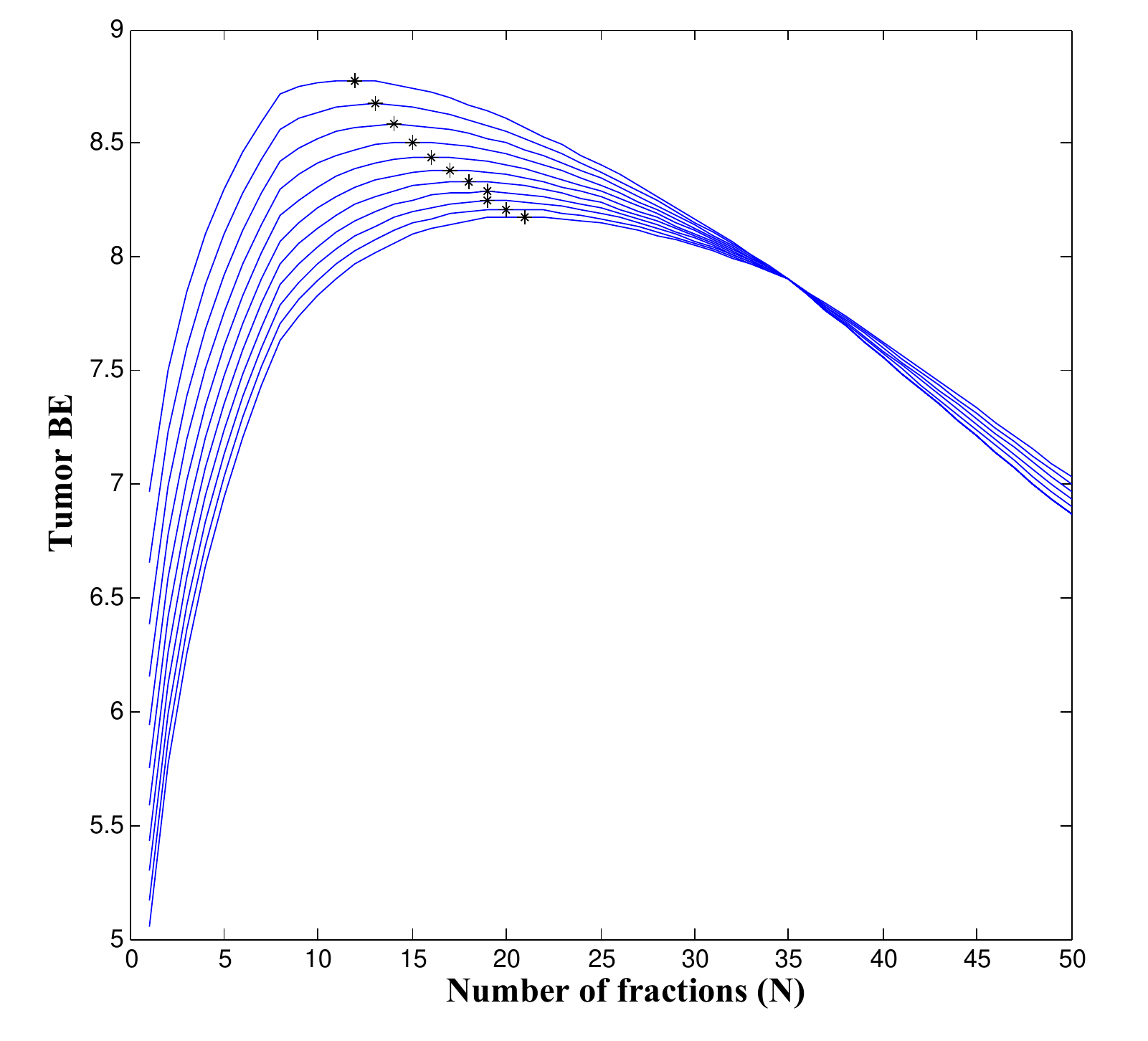}
		\caption{}
		\label{fig:7-10}
	\end{subfigure}
	\hfill
	\begin{subfigure}[b]{0.5\textwidth}
		\includegraphics[width=\textwidth]{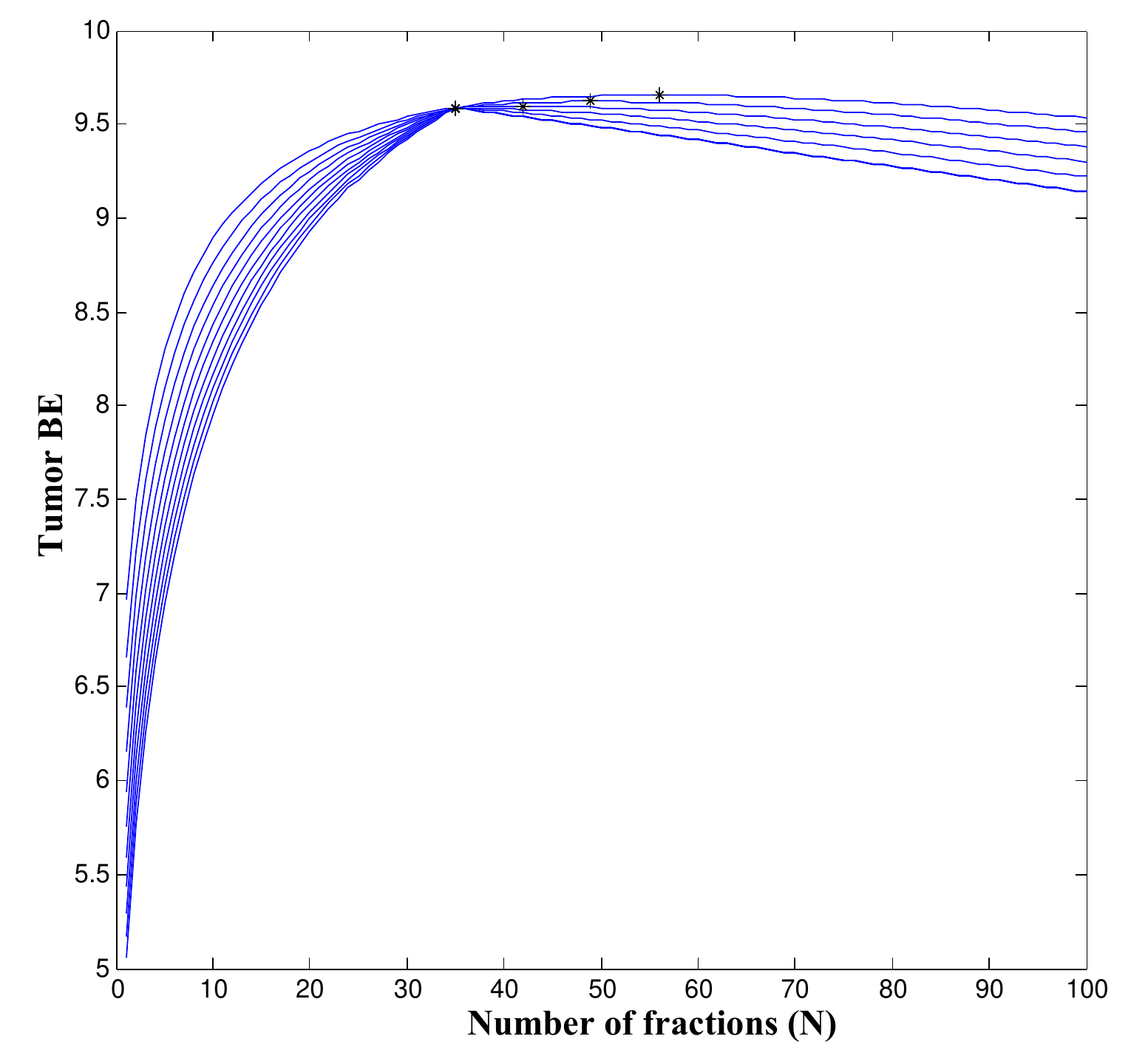}
		\caption{}
		\label{fig:7-50}
	\end{subfigure}%	
	\begin{subfigure}[b]{0.5\textwidth}
		\includegraphics[width=\textwidth]{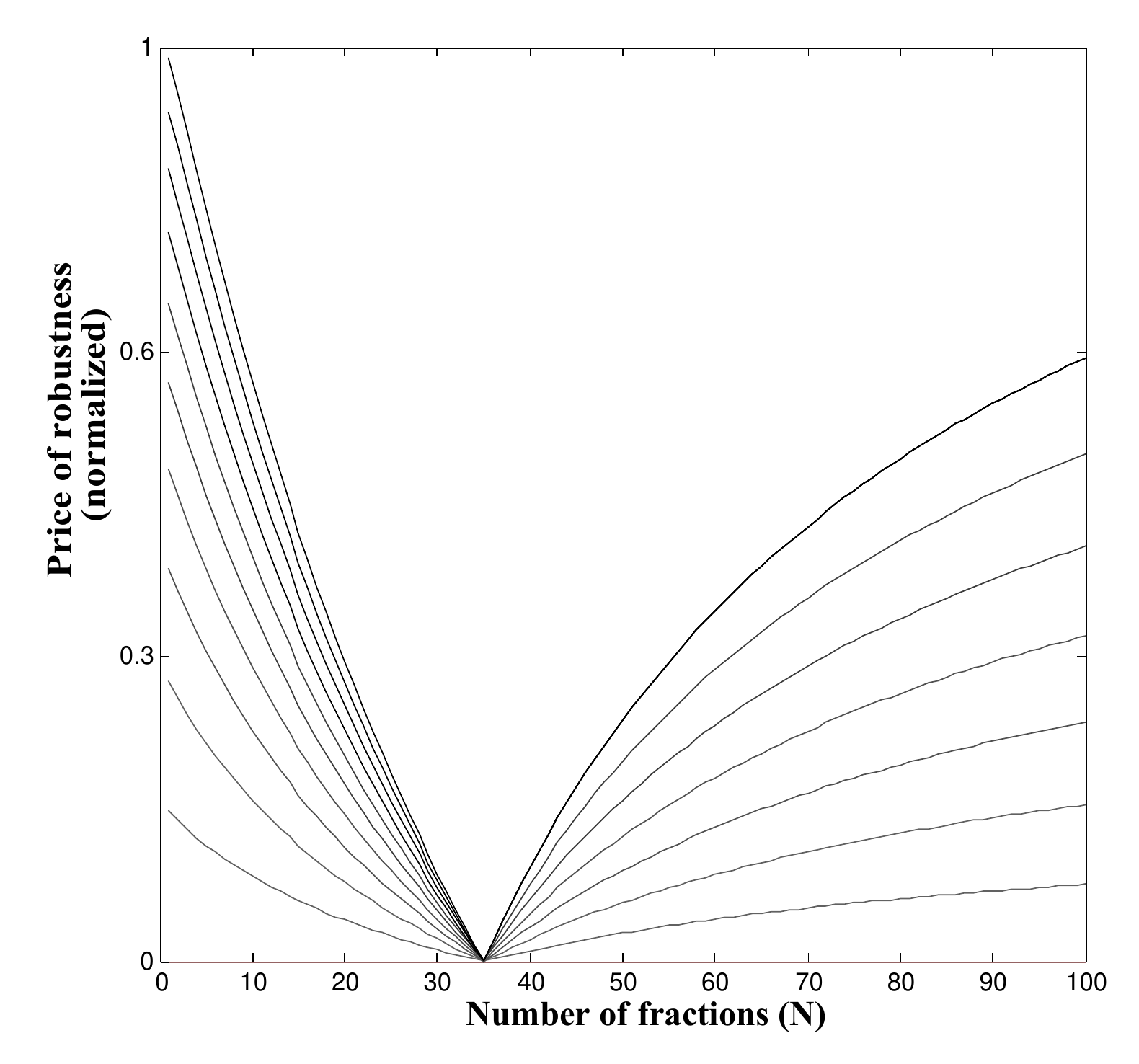}
		\caption{}
		\label{fig:PriceofRobusness}
	\end{subfigure}
	\caption{(a), (b), (c): The value of the objective as a function of $N$ for $T_{\text{lag}}=7$. The data points denoted by $\ast$ show the points $(N^*,f^*)$  in each graph. The uppermost line in each set of graphs (a), (b), and (c) represents the nominal case ($\delta=0$) and the other lines correspond to $\delta=\{0.1,0.2,\ldots,1\}$, respectively, from top to bottom. (a) $T_{\text{double}}=2$ days and $N^*=8$ for all $\delta$. (b) $T_{\text{double}}=10$ days. In this case, $N^*$ increases with increasing $\delta$. (c) $T_{\text{double}}=100$ days. In 
	this case, $N^*$ decreases with increasing $\delta$ (from top to bottom). 
	Moreover, note that in some cases (for example in (c)), some of these 
	points lie on top of each other as the optimal solution for their corresponding $\delta$ 
	are equal. (d) Normalized price of robustness as a function of $N$ when $T_{\text{lag}}=7,\ T_{\text{double}}=2$. The uppermost line represents the most uncertain case ($\delta=1$) and other lines represent $\delta=\{0.9, 0.8,\ldots,0.1\}$, respectively, from top to bottom.}
	\label{fig:objectivecomparison}
\end{figure}

For each $T_{\text{lag}}$, $\delta$ combination in Table 
\ref{table:PriceofRobustness}, the price of robustness first decreases with 
increasing $T_{\text{double}}$, reaches the 
smallest value when $T_{\text{double}}=50$ days and then increases. This trend 
is consistent 
with the corresponding trend in the difference between $N_m=35$ and $N^*$ that 
can be inferred from Table \ref{table:comparisonN&d}. Specifically, for each 
$T_{\text{lag}}$, $\delta$ combination, the magnitude of $N_m-N^*$ decreases 
with increasing
$T_{\text{double}}$, reaches about a day or two when $T_{\text{double}}=50$, 
and then increases. In fact, as we can see in Figure \ref{fig:objectivecomparison}(d), 
when $N=N_m=35$, the price of robustness is exactly zero; more strongly, we 
found that this held true irrespective of the values of $\delta, T_{\text{lag}}$, and 
$T_{\text{double}}$. A detailed algebraic proof of this fact can be 
developed, but is omitted here for brevity. Roughly speaking, the key idea in this 
proof is that when $N=N_m$, the BED constraints reduce to total dose 
constraints; this eliminates 
the dependence of the BED constraint on $\rho_m$ and hence an 
equal-dosage solution that splits the tolerance dose across $N$ fractions is 
optimal to the nominal {\emph{and}} the robust problem. Consequently, the price 
of robustness is zero. Finally, for any combination of $\delta$ and 
$T_{\text{double}}$ in Table \ref{table:PriceofRobustness}, the price of robustness 
decreases as $T_{\text{lag}}$ 
increases. Again, this is also consistent with the corresponding trend in the 
magnitude of $N_m-N^*$.

A closer look at Table \ref{table:comparisonN&d} reveals that the evolution of 
$N^*$ with $\delta$ for various fixed combinations of $T_{\text{double}}$ and
$T_{\text{lag}}$ does not exhibit a universal trend. For instance, when 
$T_{\text{lag}}=7$ days and $T_{\text{double}}=2$ days, $N^*=8$ for all $\delta$ 
(also see  
Figure \ref{fig:objectivecomparison} (a)). However, $N^*$ increases with increasing $\delta$ 
when $T_{\text{lag}}=7$ days and 
$T_{\text{double}}=10$ days (Figure \ref{fig:objectivecomparison} (b)). On the 
other hand, $N^*$ decreases as $\delta$ increases
when $T_{\text{lag}}=7$ days and $T_{\text{double}}=100$ days (Figure 
\ref{fig:objectivecomparison} (c)). Consistent with this observation, Table 
\ref{table:comparisonN&d} also shows that for each fixed combination of 
$T_{\text{double}}$ and $T_{\text{lag}}$, optimal doses do not exhibit a universal qualitative trend as a function of $\delta$.  

Finally, our robust solutions continue to exhibit qualitative trends that are well-established  
in the nominal case (see \cite{Saberian2015mathmedbio} and references therein). For instance, $N^*$ increases with increasing $T_{\text{double}}$ 
for any fixed $\delta$, $T_{\text{lag}}$ combination; similarly, $N^*$ also increases 
as $T_{\text{lag}}$ increases for any fixed $\delta$, $T_{\text{double}}$ 
combination.
\subsection{Infeasibility tests}\label{sec:infeasible}
As mentioned before, the primary motivation for the robust 
formulation is that the optimal solution obtained by solving the nominal formulation 
is guaranteed to be feasible only for nominal values of $\rho$. This means that if 
the actual $\rho$ values turn out to be different than the nominal, the nominal 
solution might become infeasible. To quantify the frequency and extent of such infeasibility, we first performed 
a set of numerical experiments where the realized values of $\rho$ were assumed 
to equal various grid-points in the uncertainty intervals around the nominal values. It turned out that while $\rho$ varied in this manner over grid-points inside the 
uncertainty interval, the nominal solution was infeasible in about 75\% of the 
cases; the robust solution of course remained 
feasible in all cases. The amount of infeasibility in 
some cases was rather large --- close to 50\%, with an average of 10.5\%. The first, second, and third quartiles were 3.99\%, 8.44\%, and 15.76\%, respectively. Given the relatively small price of 
robustness reported in Section \ref{sec:qualrobust}, this suggests that it might be 
worthwhile to implement the robust dosing schedules rather than the nominal 
ones.

In the robust counterpart, we assumed that the value of $\rho$ for each OAR 
belongs to a known interval. 
Therefore, any 
solution to our robust formulation is guaranteed to be feasible only as long 
as this assumption holds. 
Due to the uncertainty in the actual values of $\rho$, however, this assumption 
could be violated. In that case, our robust solution 
might not be feasible after all. To test the impact of this unfortunate occurrence, 
we performed numerical experiments 
where $\rho$ values were varied outside the predetermined uncertainty interval. In particular, for each uncertainty level $\delta$ and constraint $m$, five grid-points were chosen at $\tilde\rho^m=(1+\delta+\gamma)\rho^m$ and five grid-points at $\tilde \rho^m=(1-\delta-\gamma)\rho^m$, where $\gamma \in \{0.1,0.2,\ldots,0.5\}$ and $\rho^m$ denotes the nominal value. The nominal solution was infeasible in over 65\% of the cases, while the robust solution was 
infeasible in 43\% of the cases. The amount of infeasibility was found to be
statistically lower (via a pairwise t-test at the $p=0.05$ significance level) for the 
robust solution than the nominal solution over all cases. This is encouraging 
because it suggests that the robust solution might be ``more robust" than the 
nominal solution even when $\rho$ values are outside the uncertainty intervals 
(although it does not appear possible to rigorously state and prove this claim).
\subsection{Uncertainty in tumor parameters}\label{sec:tumoruncertain}
Throughout this paper, we assumed that the values of the tumor 
parameters $\alpha_0$ and $\beta_0$ were known. In this section, we 
investigate the effect of uncertainty in these tumor parameters. 
We assume that both $\tilde{\alpha}_0$ and $\tilde{\beta}_0$ belong to a known interval. That is, $\tilde{\alpha}_0 \in [\alpha_0^{\text{min}},\alpha_0^{\text{max}}]$ and $\tilde{\beta}_0 \in [\beta_0^{\text{min}},\beta_0^{\text{max}}]$. Since we are maximizing the objective function, the worst realization of the problem occurs when 
$\tilde{\alpha}_0=\alpha_0^{\text{min}}, \ \tilde{\beta}_0=\beta_0^{\text{min}}$. 
That is, the robust objective value is simply attained by 
replacing $\tilde{\alpha}_0$ and $\tilde{\beta}_0$ by their minimum values. Table \ref{table:tumoruncertainty} shows the effect of this uncertainty, assuming that $
\tilde{\alpha}_0 \in [(1-\theta)\hat{\alpha}_0,(1+\theta)\hat{\alpha}_0]$ and $
\tilde{\beta}_0 \in [(1-\theta)\hat{\beta}_0,(1+\theta)\hat{\beta}_0]$. Here, the 
nominal values $\hat{\alpha}_0$ and $\hat{\beta}_0$ were set to $0.35\ \text{Gy}$ and 
$0.035\ \text{Gy}^{-2}$, respectively, and $\theta$ was varied in the set $\{0.1, 0.2,\ldots,0.9\}$. A quick inspection of the table reveals that as the tumor uncertainty level increases, the number of fractions decreases and the dose delivered in each session increases. In other words, higher uncertainty in tumor parameters causes it to behave similar to faster-proliferating tumors.
\section{Discussion}\label{sec:discussion}
Most existing research on robust optimization in cancer radiotherapy focuses on 
incorporating uncertainty in the actual dose delivered to various anatomical 
regions of interest via intensity modulated radiation therapy (IMRT) and other 
treatment methods. Causes of this uncertainty include patient movement, say due 
to breathing, or setup errors at the time of treatment delivery (see, for instance,
\cite{bortfeldrobust2008,chan2006,mahmoudzadehbreast2015,mahmoudzadeh2015,mar2015,unkelbach2007} and references therein). 

In this paper, we provided a robust 
formulation of the fractionation problem. Perhaps more importantly, we also 
described in detail a simple method for exact solution of this robust formulation. Although our robust 
formulation is, at first glance, inevitably at least as hard as a non-convex QCQP, 
we were able to show that it can be solved to optimality by solving a few 
two-variable LPs with a few constraints each. Our numerical experiments 
provided insights into the behavior of nominal and robust dosing schedules and also quantified 
the price of robustness. %Overall, it appears that the price of robustness is small.
Overall, our comparison of the frequency and amount of infeasibility incurred by 
the nominal and the robust solutions suggests that the robust solutions are indeed 
statistically more feasible and yet pay a relatively small price of robustness. 
This could provide motivation for future investigations into the use of biological 
dose-response models such as the LQ model for planning radiation treatment as 
the uncertainty in dose-response parameters has been the main obstacle in 
widespread reliance on these models (see \cite{Hall05}).

Note that the nominal fractionation model in 
\cite{Saberian2015mathmedbio,Saberian2015} used the concept of sparing 
factors to model the doses delivered to the various OAR. For instance, if dose 
$d_t$ is given to the tumor in fraction $t$, then the dose to OAR $m\in\mathcal M$ 
equals $s_m d_t$; here, $s_m$ is a non-negative sparing factor. In this paper, we 
did not use 
such sparing factors because they would have been distracting to the main 
message of our work. We do emphasize, however, that our solution 
procedure in Section \ref{sec:solve} would work even if such sparing factors were 
included. More strongly, our solution method would work even if 
the true values of these sparing factors were unknown but were instead assumed 
to belong to a non-negative interval. 
This can be done simply by using the largest 
values of these sparing factors in our robust formulation as in \cite{badri2015v2}.

In a recent unpublished manuscript based on the doctoral dissertation of 
Saberian \cite{fatemehthesis}, Saberian et al. \cite{upaper} presented a 
spatiotemporally integrated 
formulation of the fractionation problem. The decision variables in that formulation 
were $N$ and the intensity profiles of the IMRT radiation fields employed in each 
fraction. The numbers of variables and constraints in that non-convex formulation are as large as tens of thousands. It would be interesting in the future to formulate 
the robust counterpart of that model and to devise efficient approximate solution 
methods.
\section{Acknowledgment}
Funded in part by the National Science Foundation through grant \#CMMI 1054026.

\bibliographystyle{plain}
\bibliography{robust-fractionation}

\begin{thebibliography}{10}

\bibitem{badri2015v2}
H~Badri, Y~Watanabe, and K~Leder.
\newblock Robust and probabilistic optimization of dose schedules in
  radiotherapy.
\newblock available online at \url{arxiv.org/pdf/1503.00399v2.pdf}, June 2015.

\bibitem{badri2015v1}
H~Badri, Y~Watanabe, and K~Leder.
\newblock Robust optimization of dose schedules in radiotherapy.
\newblock available online at \url{arxiv.org/pdf/1503.00399v1.pdf}, March 2015.

\bibitem{bental}
A~Ben-Tal, L~El Ghaoui, and A~Nemirovski.
\newblock {\em Robust Optimization}.
\newblock {Princeton University Press}, Princeton, NJ, USA, 2009.

\bibitem{robustreview}
D~Bertsimas, D~B Brown, and C~Caramanis.
\newblock Theory and applications of robust optimization.
\newblock {\em {SIAM Review}}, 53(3):464--501, 2011.

\bibitem{bertuzzi2013}
A~Bertuzzi, C~Bruni~F Papa, and C~Sinisgalli.
\newblock Optimal solution for a cancer radiotherapy problem.
\newblock {\em Journal of Mathematical Biology}, 66(1-2):311--349, 2013.

\bibitem{bortfeldrobust2008}
T~Bortfeld, T~C~Y Chan, A~Trofimov, and J~N Tsitsiklis.
\newblock Robust management of motion uncertainty in intensity modulated
  radiation therapy.
\newblock {\em Operations Research}, 56(6):1461--1473, 2008.

\bibitem{bortfeld13}
T~Bortfeld, J~Ramakrishnan, J~N Tsitsiklis, and J~Unkelbach.
\newblock Optimization of radiotherapy fractionation schedules in the presence
  of tumor repopluation.
\newblock forthcoming in {INFORMS Journal on Computing}, prepring available at
  \url{http://pages.discovery.wisc.edu/~jramakrishnan/BRT2015_repop.pdf}, May
  2015.

\bibitem{chan2006}
T~C~Y Chan, T~Bortfeld, and J~Tsitsiklis.
\newblock A robust approach to {IMRT} optimization.
\newblock {\em Physics in Medicine and Biology}, 51:2567--2583, 2006.

\bibitem{fowler1990}
J~F Fowler.
\newblock How worthwhile are short schedules in radiotherapy?: A series of
  exploratory calculations.
\newblock {\em Radiotherapy and Oncology}, 18(2):165--181, 1990.

\bibitem{fowler01}
J~F Fowler.
\newblock Biological factors influencing optimum fractionation in radiation
  therapy.
\newblock {\em Acta Oncologica}, 40(6):712--717, 2001.

\bibitem{fowler2007}
J~F Fowler.
\newblock Is there an optimal overall time for head and neck radiotherapy? a
  review with new modeling.
\newblock {\em Clinical Oncology}, 19(1):8--27, 2007.

\bibitem{fowler2008}
J~F Fowler.
\newblock {Optimum overall times II: Extended modelling for head and neck
  radiotherapy}.
\newblock {\em Clinical Oncology}, 20(2):113--126, 2008.

\bibitem{Fowler1995}
J~F Fowler and M~A Ritter.
\newblock A rationale for fractionation for slowly proliferating tumors such as
  prostatic adenocarcinoma.
\newblock {\em International Journal of Radiation Oncology Biology Physics},
  32(2):521--529, 1995.

\bibitem{Hall05}
E~J Hall and A~J Giaccia.
\newblock {\em Radiobiology for the {R}adiologist}.
\newblock Lippincott Williams \& Wilkins, Philadelphia, Pennsylvania, USA,
  2005.

\bibitem{Jones1995}
B~Jones, L~T Tan, and R~G Dale.
\newblock Derivation of the optimum dose per fraction from the linear quadratic
  model.
\newblock {\em The British Journal of Radiology}, 68(812):894--902, 1995.

\bibitem{mahmoudzadehbreast2015}
H~Mahmoudzadeh, J~Lee, T~C~Y Chan, and T~G Purdie.
\newblock {Robust optimization methods for cardiac sparing in tangential breast
  IMRT}.
\newblock {\em Medical Physics}, 42(5):2212, 2015.

\bibitem{mahmoudzadeh2015}
H~Mahmoudzadeh, T~G Purdie, and T~C~Y Chan.
\newblock Constraint generation methods for robust optimization in radiation
  therapy.
\newblock forthcoming in Operations Research for Health Care, June 2015.

\bibitem{mar2015}
P~A Mar and T~C~Y Chan.
\newblock Adaptive and robust radiation therapy in the presence of drift.
\newblock {\em Physics in Medicine and Biology}, 60(9):3599--3615, 2015.

\bibitem{QUANTEC}
L~B Marks, E~D Yorke, A~Jackson, R~K~Ten Haken, L~S Constine, A~Eisbruch, S~M
  Bentzen, J~Nam, and J~O Deasy.
\newblock Use of normal tissue complication probability models in the clinic.
\newblock {\em International Journal of Radiation Oncology Biology Physics},
  76(3):S10--S19, 2010.

\bibitem{Mizuta2012}
M~Mizuta, S~Takao, H~Date, N~Kishimoto, K~L Sutherland, R~Onimaru, and
  H~Shirato.
\newblock A mathematical study to select fractionation regimen based on
  physical dose distribution and the linear-quadratic model.
\newblock {\em International Journal of Radiation Oncology Biology Physics},
  84(3):829 -- 833, 2012.

\bibitem{qi2012}
X~S Qi, Q~Yang, S~P Lee, X~A Li, and D~Wang.
\newblock An estimation of radiobiological parameters for head-and-neck cancer
  cells and the clinical implications.
\newblock {\em Cancers}, 4:566--580, 2012.

\bibitem{Rockwell:1998ve}
S~Rockwell.
\newblock Experimental radiotherapy: a brief history.
\newblock {\em Radiation Research}, 150(Supplement):S157--S169, November 1998.

\bibitem{fatemehthesis}
F~Saberian.
\newblock {\em Convex and Dynamic Optimization with Learning for Adaptive
  Biologically Conformal Radiotherapy}.
\newblock unpublished doctoral dissertation, University of Washington,
  Industrial and Systems Engineering, May 2015.

\bibitem{Saberian2015mathmedbio}
F~Saberian, A~Ghate, and M~Kim.
\newblock Optimal fractionation in radiotherapy with multiple normal tissues.
\newblock forthcoming in Mathematical Medicine and Biology, online preprint
  available at doi: 10.1093/imammb/dqv015, May 2015.

\bibitem{upaper}
F~Saberian, A~Ghate, and M~Kim.
\newblock Spatiotemporally integrated fractionation in radiotherapy.
\newblock under review at {INFORMS Journal on Computing}, preprint available at
  \url{http://faculty.washington.edu/archis/upaper-apr-2015.pdf}, April 2015.

\bibitem{Saberian2015}
F~Saberian, A~Ghate, and M~Kim.
\newblock A two-variable linear program solves the standard linear--quadratic
  formulation of the fractionation problem in cancer radiotherapy.
\newblock {\em {Operations Research Letters}}, 43(3):254 -- 258, 2015.

\bibitem{unkelbach2007}
J~Unkelbach, T~C~Y Chan, and T~Bortfeld.
\newblock Accounting for range uncertainties in the optimization of intensity
  modulated proton therapy.
\newblock {\em Physics in Medicine and Biology}, 52:2755--2773, 2007.

\bibitem{unkelbach13}
J~Unkelbach, D~Craft, E~Saleri, J~Ramakrishnan, and T~Bortfeld.
\newblock The dependence of optimal fractionation schemes on the spatial dose
  distribution.
\newblock {\em Physics in Medicine and Biology}, 58(1):159--167, 2013.

\bibitem{Yang2005fractionation}
Y~Yang and L~Xing.
\newblock Optimization of radiotherapy dose-time fractionation with
  consideration of tumor specific biology.
\newblock {\em Medical Physics}, 32(12):3666--3677, 2005.

\end{thebibliography}

\newpage
\begin{longtable}{|l|l|l|l|l|l|l|l|l|}
	\hline
	$T_{\text{lag}}=7$ & \multicolumn{8}{c|}{$T_{\text{double}}$} \\ 
	\hline
	$\delta$ & 2 & 8 & 10 & 20 & 40 & 50 & 80 & 100 \\ 
	\hline
	0.1 & 1.74\% & 1.44\% & 1.21\% & 0.56\% & 0.04\% & 0.01\% & 0.25\% & 0.36\% \\ 
	0.2 & 3.34\% & 2.67\% & 2.23\% & 0.98\% & 0.04\% & 0.01\% & 0.39\% & 0.62\% \\ 
	0.3 & 4.80\% & 3.74\% & 3.10\% & 1.31\% & 0.04\% & 0.01\% & 0.39\% & 0.70\% \\ 
	0.4 & 6.14\% & 4.66\% & 3.85\% & 1.56\% & 0.04\% & 0.01\% & 0.39\% & 0.70\% \\ 
	0.5 & 7.38\% & 5.48\% & 4.51\% & 1.74\% & 0.04\% & 0.01\% & 0.39\% & 0.70\% \\ 
	0.6 & 8.53\% & 6.22\% & 5.09\% & 1.87\% & 0.04\% & 0.01\% & 0.39\% & 0.70\% \\ 
	0.7 & 9.61\% & 6.88\% & 5.61\% & 1.95\% & 0.04\% & 0.01\% & 0.39\% & 0.70\% \\ 
	0.8 & 10.61\% & 7.47\% & 6.07\% & 2.00\% & 0.04\% & 0.01\% & 0.39\% & 0.70\% \\ 
	0.9 & 11.54\% & 8.01\% & 6.47\% & 2.02\% & 0.04\% & 0.01\% & 0.39\% & 0.70\% \\ 
	1 & 12.42\% & 8.50\% & 6.84\% & 2.02\% & 0.04\% & 0.01\% & 0.39\% & 0.70\% \\ 
	\hline
	\hline
	$T_{\text{lag}}=14$ & \multicolumn{8}{c|}{$T_{\text{double}}$} \\ 
	\hline
	$\delta$ & 2 & 8 & 10 & 20 & 40 & 50 & 80 & 100 \\ 
	\hline
	0.1 & 0.92\% & 0.92\% & 0.92\% & 0.54\% & 0.04\% & 0.01\% & 0.25\% & 0.36\% \\ 
	0.2 & 1.78\% & 1.78\% & 1.78\% & 0.96\% & 0.04\% & 0.01\% & 0.38\% & 0.61\% \\ 
	0.3 & 2.59\% & 2.59\% & 2.59\% & 1.28\% & 0.04\% & 0.01\% & 0.39\% & 0.70\% \\ 
	0.4 & 3.34\% & 3.34\% & 3.30\% & 1.52\% & 0.04\% & 0.01\% & 0.39\% & 0.70\% \\ 
	0.5 & 4.05\% & 4.05\% & 3.93\% & 1.69\% & 0.04\% & 0.01\% & 0.39\% & 0.70\% \\ 
	0.6 & 4.71\% & 4.71\% & 4.48\% & 1.82\% & 0.04\% & 0.01\% & 0.39\% & 0.70\% \\ 
	0.7 & 5.34\% & 5.34\% & 4.97\% & 1.90\% & 0.04\% & 0.01\% & 0.39\% & 0.70\% \\ 
	0.8 & 5.93\% & 5.89\% & 5.41\% & 1.95\% & 0.04\% & 0.01\% & 0.39\% & 0.70\% \\ 
	0.9 & 6.49\% & 6.41\% & 5.80\% & 1.96\% & 0.04\% & 0.01\% & 0.39\% & 0.70\% \\ 
	1 & 7.02\% & 6.87\% & 6.14\% & 1.96\% & 0.04\% & 0.01\% & 0.39\% & 0.70\% \\ 
	\hline
	\hline
	$T_{\text{lag}}=21$ & \multicolumn{8}{c|}{$T_{\text{double}}$} \\ 
	\hline
	$\delta$ & 2 & 8 & 10 & 20 & 40 & 50 & 80 & 100 \\ 
	\hline
	0.1 & 0.47\% & 0.47\% & 0.47\% & 0.47\% & 0.04\% & 0.01\% & 0.25\% & 0.36\% \\ 
	0.2 & 0.92\% & 0.92\% & 0.92\% & 0.88\% & 0.04\% & 0.01\% & 0.38\% & 0.61\% \\ 
	0.3 & 1.34\% & 1.34\% & 1.34\% & 1.19\% & 0.04\% & 0.01\% & 0.39\% & 0.69\% \\ 
	0.4 & 1.74\% & 1.74\% & 1.74\% & 1.42\% & 0.04\% & 0.01\% & 0.39\% & 0.69\% \\ 
	0.5 & 2.12\% & 2.12\% & 2.12\% & 1.60\% & 0.04\% & 0.01\% & 0.39\% & 0.69\% \\ 
	0.6 & 2.48\% & 2.48\% & 2.48\% & 1.72\% & 0.04\% & 0.01\% & 0.39\% & 0.69\% \\ 
	0.7 & 2.82\% & 2.82\% & 2.82\% & 1.80\% & 0.04\% & 0.01\% & 0.39\% & 0.69\% \\ 
	0.8 & 3.14\% & 3.14\% & 3.14\% & 1.85\% & 0.04\% & 0.01\% & 0.39\% & 0.69\% \\ 
	0.9 & 3.45\% & 3.45\% & 3.45\% & 1.86\% & 0.04\% & 0.01\% & 0.39\% & 0.69\% \\ 
	1 & 3.75\% & 3.75\% & 3.75\% & 1.86\% & 0.04\% & 0.01\% & 0.39\% & 0.69\% \\ 
	\hline
	\hline
	$T_{\text{lag}}=28$ & \multicolumn{8}{c|}{$T_{\text{double}}$} \\ 
	\hline
	$\delta$ & 2 & 8 & 10 & 20 & 40 & 50 & 80 & 100 \\ 
	\hline
	0.1 & 0.18\% & 0.18\% & 0.18\% & 0.18\% & 0.04\% & 0.01\% & 0.25\% & 0.36\% \\ 
	0.2 & 0.35\% & 0.35\% & 0.35\% & 0.35\% & 0.04\% & 0.01\% & 0.38\% & 0.61\% \\ 
	0.3 & 0.52\% & 0.52\% & 0.52\% & 0.52\% & 0.04\% & 0.01\% & 0.38\% & 0.69\% \\ 
	0.4 & 0.67\% & 0.67\% & 0.67\% & 0.67\% & 0.04\% & 0.01\% & 0.38\% & 0.69\% \\ 
	0.5 & 0.82\% & 0.82\% & 0.82\% & 0.82\% & 0.04\% & 0.01\% & 0.38\% & 0.69\% \\ 
	0.6 & 0.97\% & 0.97\% & 0.97\% & 0.94\% & 0.04\% & 0.01\% & 0.38\% & 0.69\% \\ 
	0.7 & 1.10\% & 1.10\% & 1.10\% & 1.02\% & 0.04\% & 0.01\% & 0.38\% & 0.69\% \\ 
	0.8 & 1.23\% & 1.23\% & 1.23\% & 1.07\% & 0.04\% & 0.01\% & 0.38\% & 0.69\% \\ 
	0.9 & 1.36\% & 1.36\% & 1.36\% & 1.08\% & 0.04\% & 0.01\% & 0.38\% & 0.69\% \\ 
	1 & 1.48\% & 1.48\% & 1.48\% & 1.08\% & 0.04\% & 0.01\% & 0.38\% & 0.69\% \\ 
	\hline
	\hline
	$T_{\text{lag}}=35$ & \multicolumn{8}{c|}{$T_{\text{double}}$} \\ 
	\hline
	$\delta$ & 2 & 8 & 10 & 20 & 40 & 50 & 80 & 100 \\ 
	\hline
	0.1 & 0.03\% & 0.03\% & 0.03\% & 0.03\% & 0.03\% & 0.03\% & 0.25\% & 0.36\% \\ 
	0.2 & 0.06\% & 0.06\% & 0.06\% & 0.06\% & 0.06\% & 0.06\% & 0.38\% & 0.61\% \\ 
	0.3 & 0.08\% & 0.08\% & 0.08\% & 0.08\% & 0.08\% & 0.09\% & 0.41\% & 0.69\% \\ 
	0.4 & 0.12\% & 0.12\% & 0.12\% & 0.12\% & 0.12\% & 0.12\% & 0.44\% & 0.72\% \\ 
	0.5 & 0.15\% & 0.15\% & 0.15\% & 0.15\% & 0.15\% & 0.15\% & 0.47\% & 0.76\% \\ 
	0.6 & 0.15\% & 0.15\% & 0.15\% & 0.15\% & 0.15\% & 0.15\% & 0.47\% & 0.76\% \\ 
	0.7 & 0.15\% & 0.15\% & 0.15\% & 0.15\% & 0.15\% & 0.15\% & 0.47\% & 0.76\% \\ 
	0.8 & 0.15\% & 0.15\% & 0.15\% & 0.15\% & 0.15\% & 0.15\% & 0.47\% & 0.76\% \\ 
	0.9 & 0.15\% & 0.15\% & 0.15\% & 0.15\% & 0.15\% & 0.15\% & 0.47\% & 0.76\% \\ 
	1 & 0.15\% & 0.15\% & 0.15\% & 0.15\% & 0.15\% & 0.15\% & 0.47\% & 0.76\% \\ 
	\hline
	\hline
	\caption{The price of robustness for different combinations of $T_{\text{lag}}$, $T_{\text{double}}$, and $\delta$. The price of robustness equals 
	$(\frac{g^*-f^*}{g^*})\times 100\%$, where $f^*$ and $g^*$ are the optimal values of the robust and the nominal formulations, respectively.}
	\label{table:PriceofRobustness}
\end{longtable}

\begin{landscape}
	\small
	\begin{longtable}{|l|l|l|l|l|l|l|l|l|}
		\hline
		$T_{\text{lag}}=7$ & \multicolumn{8}{c|}{$T_{\text{double}}$} \\ 
		\hline
		$\delta$ & 2 & 8 & 10 & 20 & 40 & 50 & 80 & 100 \\ 
		\hline
		0 & (2.49, 8) & (2.10,10) & (1.82,12) & (1.20,20) & (0.80,32) & (0.71,37) & (0.55,49) & (0.49,56) \\ 
		0.1 & (2.46, 8) & (1.93,11) & (1.69,13) & (1.11,22) & (0.74,35) & (0.74,35) & (0.62,43) & (0.55,49) \\ 
		0.2 & (2.42, 8) & (1.79,12) & (1.59,14) & (1.02,24) & (0.74,35) & (0.74,35) & (0.71,37) & (0.63,42) \\ 
		0.3 & (2.39, 8) & (1.67,13) & (1.49,15) & (0.96,26) & (0.74,35) & (0.74,35) & (0.74,35) & (0.74,35) \\ 
		0.4 & (2.36, 8) & (1.65,13) & (1.41,16) & (0.92,27) & (0.74,35) & (0.74,35) & (0.74,35) & (0.74,35) \\ 
		0.5 & (2.34, 8) & (1.55,14) & (1.34,17) & (0.87,29) & (0.74,35) & (0.74,35) & (0.74,35) & (0.74,35) \\ 
		0.6 & (2.31, 8) & (1.46,15) & (1.27,18) & (0.82,31) & (0.74,35) & (0.74,35) & (0.74,35) & (0.74,35) \\ 
		0.7 & (2.29, 8) & (1.45,15) & (1.21,19) & (0.80,32) & (0.74,35) & (0.74,35) & (0.74,35) & (0.74,35) \\ 
		0.8 & (2.27, 8) & (1.38,16) & (1.21,19) & (0.76,34) & (0.74,35) & (0.74,35) & (0.74,35) & (0.74,35) \\ 
		0.9 & (2.25, 8) & (1.31,17) & (1.16,20) & (0.74,35) & (0.74,35) & (0.74,35) & (0.74,35) & (0.74,35) \\ 
		1 & (2.23, 8) & (1.30,17) & (1.11,21) & (0.74,35) & (0.74,35) & (0.74,35) & (0.74,35) & (0.74,35) \\ 
		\hline
		\hline
		$T_{\text{lag}}=14$ & \multicolumn{8}{c|}{$T_{\text{double}}$} \\ 
		\hline
		$\delta$ & 2 & 8 & 10 & 20 & 40 & 50 & 80 & 100 \\ 
		\hline
		0 & (1.53,15) & (1.53,15) & (1.53,15) & (1.20,20) & (0.80,32) & (0.71,37) & (0.55,49) & (0.49,56) \\ 
		0.1 & (1.51,15) & (1.51,15) & (1.51,15) & (1.11,22) & (0.74,35) & (0.74,35) & (0.62,43) & (0.55,49) \\ 
		0.2 & (1.50,15) & (1.50,15) & (1.50,15) & (1.02,24) & (0.74,35) & (0.74,35) & (0.71,37) & (0.63,42) \\ 
		0.3 & (1.49,15) & (1.49,15) & (1.49,15) & (0.96,26) & (0.74,35) & (0.74,35) & (0.74,35) & (0.74,35) \\ 
		0.4 & (1.48,15) & (1.48,15) & (1.41,16) & (0.92,27) & (0.74,35) & (0.74,35) & (0.74,35) & (0.74,35) \\ 
		0.5 & (1.47,15) & (1.47,15) & (1.34,17) & (0.87,29) & (0.74,35) & (0.74,35) & (0.74,35) & (0.74,35) \\ 
		0.6 & (1.46,15) & (1.46,15) & (1.27,18) & (0.82,31) & (0.74,35) & (0.74,35) & (0.74,35) & (0.74,35) \\ 
		0.7 & (1.45,15) & (1.45,15) & (1.21,19) & (0.80,32) & (0.74,35) & (0.74,35) & (0.74,35) & (0.74,35) \\ 
		0.8 & (1.45,15) & (1.38,16) & (1.21,19) & (0.76,34) & (0.74,35) & (0.74,35) & (0.74,35) & (0.74,35) \\ 
		0.9 & (1.44,15) & (1.31,17) & (1.16,20) & (0.74,35) & (0.74,35) & (0.74,35) & (0.74,35) & (0.74,35) \\ 
		1 & (1.43,15) & (1.30,17) & (1.11,21) & (0.74,35) & (0.74,35) & (0.74,35) & (0.74,35) & (0.74,35) \\ 
		\hline
		\hline
		$T_{\text{lag}}=21$ & \multicolumn{8}{c|}{$T_{\text{double}}$} \\ 
		\hline
		$\delta$ & 2 & 8 & 10 & 20 & 40 & 50 & 80 & 100 \\ 
		\hline
		0 & (1.11,22) & (1.11,22) & (1.11,22) & (1.11,22) & (0.80,32) & (0.71,37) & (0.55,49) & (0.49,56) \\ 
		0.1 & (1.11,22) & (1.11,22) & (1.11,22) & (1.11,22) & (0.74,35) & (0.74,35) & (0.62,43) & (0.55,49) \\ 
		0.2 & (1.10,22) & (1.10,22) & (1.10,22) & (1.02,24) & (0.74,35) & (0.74,35) & (0.71,37) & (0.63,42) \\ 
		0.3 & (1.10,22) & (1.10,22) & (1.10,22) & (0.96,26) & (0.74,35) & (0.74,35) & (0.74,35) & (0.74,35) \\ 
		0.4 & (1.09,22) & (1.09,22) & (1.09,22) & (0.92,27) & (0.74,35) & (0.74,35) & (0.74,35) & (0.74,35) \\ 
		0.5 & (1.09,22) & (1.09,22) & (1.09,22) & (0.87,29) & (0.74,35) & (0.74,35) & (0.74,35) & (0.74,35) \\ 
		0.6 & (1.09,22) & (1.09,22) & (1.09,22) & (0.82,31) & (0.74,35) & (0.74,35) & (0.74,35) & (0.74,35) \\ 
		0.7 & (1.08,22) & (1.08,22) & (1.08,22) & (0.80,32) & (0.74,35) & (0.74,35) & (0.74,35) & (0.74,35) \\ 
		0.8 & (1.08,22) & (1.08,22) & (1.08,22) & (0.76,34) & (0.74,35) & (0.74,35) & (0.74,35) & (0.74,35) \\ 
		0.9 & (1.08,22) & (1.08,22) & (1.08,22) & (0.74,35) & (0.74,35) & (0.74,35) & (0.74,35) & (0.74,35) \\ 
		1 & (1.07,22) & (1.07,22) & (1.07,22) & (0.74,35) & (0.74,35) & (0.74,35) & (0.74,35) & (0.74,35) \\ 
		\hline
		\hline
		$T_{\text{lag}}=28$ & \multicolumn{8}{c|}{$T_{\text{double}}$} \\ 
		\hline
		$\delta$ & 2 & 8 & 10 & 20 & 40 & 50 & 80 & 100 \\ 
		\hline
		0 & (0.88,29) & (0.88,29) & (0.88,29) & (0.88,29) & (0.80,32) & (0.71,37) & (0.55,49) & (0.49,56) \\ 
		0.1 & (0.87,29) & (0.87,29) & (0.87,29) & (0.87,29) & (0.74,35) & (0.74,35) & (0.62,43) & (0.55,49) \\ 
		0.2 & (0.87,29) & (0.87,29) & (0.87,29) & (0.87,29) & (0.74,35) & (0.74,35) & (0.71,37) & (0.63,42) \\ 
		0.3 & (0.87,29) & (0.87,29) & (0.87,29) & (0.87,29) & (0.74,35) & (0.74,35) & (0.74,35) & (0.74,35) \\ 
		0.4 & (0.87,29) & (0.87,29) & (0.87,29) & (0.87,29) & (0.74,35) & (0.74,35) & (0.74,35) & (0.74,35) \\ 
		0.5 & (0.87,29) & (0.87,29) & (0.87,29) & (0.87,29) & (0.74,35) & (0.74,35) & (0.74,35) & (0.74,35) \\ 
		0.6 & (0.87,29) & (0.87,29) & (0.87,29) & (0.82,31) & (0.74,35) & (0.74,35) & (0.74,35) & (0.74,35) \\ 
		0.7 & (0.87,29) & (0.87,29) & (0.87,29) & (0.80,32) & (0.74,35) & (0.74,35) & (0.74,35) & (0.74,35) \\ 
		0.8 & (0.87,29) & (0.87,29) & (0.87,29) & (0.76,34) & (0.74,35) & (0.74,35) & (0.74,35) & (0.74,35) \\ 
		0.9 & (0.87,29) & (0.87,29) & (0.87,29) & (0.74,35) & (0.74,35) & (0.74,35) & (0.74,35) & (0.74,35) \\ 
		1 & (0.86,29) & (0.86,29) & (0.86,29) & (0.74,35) & (0.74,35) & (0.74,35) & (0.74,35) & (0.74,35) \\ 
		\hline
		\hline
		$T_{\text{lag}}=35$ & \multicolumn{8}{c|}{$T_{\text{double}}$} \\ 
		\hline
		$\delta$ & 2 & 8 & 10 & 20 & 40 & 50 & 80 & 100 \\ 
		\hline
		0 & (0.72,36) & (0.72,36) & (0.72,36) & (0.72,36) & (0.72,36) & (0.71,37) & (0.55,49) & (0.49,56) \\ 
		0.1 & (0.72,36) & (0.72,36) & (0.72,36) & (0.72,36) & (0.72,36) & (0.72,36) & (0.62,43) & (0.55,49) \\ 
		0.2 & (0.72,36) & (0.72,36) & (0.72,36) & (0.72,36) & (0.72,36) & (0.72,36) & (0.71,37) & (0.63,42) \\ 
		0.3 & (0.72,36) & (0.72,36) & (0.72,36) & (0.72,36) & (0.72,36) & (0.72,36) & (0.72,36) & (0.72,36) \\ 
		0.4 & (0.72,36) & (0.72,36) & (0.72,36) & (0.72,36) & (0.72,36) & (0.72,36) & (0.72,36) & (0.72,36) \\ 
		0.5 & (0.74,35) & (0.74,35) & (0.74,35) & (0.74,35) & (0.74,35) & (0.74,35) & (0.74,35) & (0.74,35) \\ 
		0.6 & (1.44,0.70,36)* & (1.44,0.70,36)* & (1.44,0.70,36)* & (1.44,0.70,36)* & (1.44,0.70,36)* & (1.44,0.70,36)* & (1.44,0.70,36)* & (1.44,0.70,36)* \\ 
		0.7 & (1.44,0.70,36)* & (1.44,0.70,36)* & (1.44,0.70,36)* & (1.44,0.70,36)* & (1.44,0.70,36)* & (1.44,0.70,36)* & (1.44,0.70,36)* & (1.44,0.70,36)* \\ 
		0.8 & (0.74,35) & (0.74,35) & (0.74,35) & (0.74,35) & (0.74,35) & (0.74,35) & (0.74,35) & (0.74,35) \\ 
		0.9 & (1.44,0.70,36)* & (1.44,0.70,36)* & (1.44,0.70,36)* & (1.44,0.70,36)* & (1.44,0.70,36)* & (1.44,0.70,36)* & (1.44,0.70,36)* & (1.44,0.70,36)* \\ 
		1 & (1.44,0.70,36)* & (1.44,0.70,36)* & (1.44,0.70,36)* & (1.44,0.70,36)* & (1.44,0.70,36)* & (1.44,0.70,36)* & (1.44,0.70,36)* & (1.44,0.70,36)* \\ 
		\hline
		\hline
		\caption{The optimal solution ($d^*,N^*$) of the robust and nominal models for different combinations of $T_{\text{lag}}, T_{\text{double}}$, and 
		$\delta$. The first row in each section of the table shows the solution of the nominal case ($\delta=0$). The cases marked with an asterisk 
		($\ast$) yield unequal-dosage solutions characterized by two dose values $(q,p)$ (recall formulas (\ref{eqn:p})-(\ref{eqn:q})) and an $N^*$ value in that order.} 
		\label{table:comparisonN&d}
	\end{longtable}
\end{landscape}

\begin{landscape}
	\small
\begin{longtable}{|l|l|l|l|l|l|l|l|l|l|l|l|}

			\hline
			\multicolumn{10}{|c|}{$T_{double}=8$} \\ 
			\hline
			$T_{lag}=7$ & \multicolumn{9}{c|}{$\theta$} \\ 
			\hline
			$\delta$ & 0.1 & 0.2 & 0.3 & 0.4 & 0.5 & 0.6 & 0.7 & 0.8 & 0.9 \\ 
			0 & (2.28, 9) & (2.49, 8) & (2.49, 8) & (2.49, 8) & (2.49, 8) & (2.49, 8) & (2.49, 8) & (2.49, 8) & (2.49, 8) \\ 
			0.1 & (2.08,10) & (2.25, 9) & (2.46, 8) & (2.46, 8) & (2.46, 8) & (2.46, 8) & (2.46, 8) & (2.46, 8) & (2.46, 8) \\ 
			0.2 & (1.91,11) & (2.05,10) & (2.22, 9) & (2.42, 8) & (2.42, 8) & (2.42, 8) & (2.42, 8) & (2.42, 8) & (2.42, 8) \\ 
			0.3 & (1.89,11) & (2.03,10) & (2.19, 9) & (2.39, 8) & (2.39, 8) & (2.39, 8) & (2.39, 8) & (2.39, 8) & (2.39, 8) \\ 
			0.4 & (1.75,12) & (1.87,11) & (2.01,10) & (2.36, 8) & (2.36, 8) & (2.36, 8) & (2.36, 8) & (2.36, 8) & (2.36, 8) \\ 
			0.5 & (1.64,13) & (1.74,12) & (1.99,10) & (2.15, 9) & (2.34, 8) & (2.34, 8) & (2.34, 8) & (2.34, 8) & (2.34, 8) \\ 
			0.6 & (1.63,13) & (1.73,12) & (1.84,11) & (2.13, 9) & (2.31, 8) & (2.31, 8) & (2.31, 8) & (2.31, 8) & (2.31, 8) \\ 
			0.7 & (1.53,14) & (1.62,13) & (1.83,11) & (1.96,10) & (2.29, 8) & (2.29, 8) & (2.29, 8) & (2.29, 8) & (2.29, 8) \\ 
			0.8 & (1.45,15) & (1.61,13) & (1.70,12) & (1.94,10) & (2.27, 8) & (2.27, 8) & (2.27, 8) & (2.27, 8) & (2.27, 8) \\ 
			0.9 & (1.44,15) & (1.51,14) & (1.69,12) & (1.93,10) & (2.25, 8) & (2.25, 8) & (2.25, 8) & (2.25, 8) & (2.25, 8) \\ 
			1 & (1.36,16) & (1.50,14) & (1.68,12) & (1.79,11) & (2.06, 9) & (2.23, 8) & (2.23, 8) & (2.23, 8) & (2.23, 8) \\ 
			\hline
			\multicolumn{10}{|c|}{$T_{double}=10$} \\ 
			\hline
			$T_{lag}=7$ & \multicolumn{9}{c|}{$\theta$} \\ 
			\hline
			$\delta$ & 0.1 & 0.2 & 0.3 & 0.4 & 0.5 & 0.6 & 0.7 & 0.8 & 0.9 \\ 
			0 & (1.95,11) & (2.10,10) & (2.28, 9) & (2.49, 8) & (2.49, 8) & (2.49, 8) & (2.49, 8) & (2.49, 8) & (2.49, 8) \\ 
			0.1 & (1.80,12) & (1.93,11) & (2.08,10) & (2.46, 8) & (2.46, 8) & (2.46, 8) & (2.46, 8) & (2.46, 8) & (2.46, 8) \\ 
			0.2 & (1.68,13) & (1.79,12) & (2.05,10) & (2.22, 9) & (2.42, 8) & (2.42, 8) & (2.42, 8) & (2.42, 8) & (2.42, 8) \\ 
			0.3 & (1.57,14) & (1.67,13) & (1.89,11) & (2.03,10) & (2.39, 8) & (2.39, 8) & (2.39, 8) & (2.39, 8) & (2.39, 8) \\ 
			0.4 & (1.48,15) & (1.65,13) & (1.75,12) & (2.01,10) & (2.17, 9) & (2.36, 8) & (2.36, 8) & (2.36, 8) & (2.36, 8) \\ 
			0.5 & (1.40,16) & (1.55,14) & (1.64,13) & (1.86,11) & (2.15, 9) & (2.34, 8) & (2.34, 8) & (2.34, 8) & (2.34, 8) \\ 
			0.6 & (1.39,16) & (1.46,15) & (1.63,13) & (1.84,11) & (1.97,10) & (2.31, 8) & (2.31, 8) & (2.31, 8) & (2.31, 8) \\ 
			0.7 & (1.32,17) & (1.45,15) & (1.53,14) & (1.71,12) & (1.96,10) & (2.29, 8) & (2.29, 8) & (2.29, 8) & (2.29, 8) \\ 
			0.8 & (1.26,18) & (1.38,16) & (1.52,14) & (1.70,12) & (1.94,10) & (2.27, 8) & (2.27, 8) & (2.27, 8) & (2.27, 8) \\ 
			0.9 & (1.25,18) & (1.31,17) & (1.44,15) & (1.60,13) & (1.80,11) & (2.25, 8) & (2.25, 8) & (2.25, 8) & (2.25, 8) \\ 
			1 & (1.20,19) & (1.30,17) & (1.43,15) & (1.59,13) & (1.79,11) & (2.06, 9) & (2.23, 8) & (2.23, 8) & (2.23, 8) \\ 
			\hline
			\multicolumn{10}{|c|}{$T_{double}=20$} \\ 
			\hline
			$T_{lag}=7$ & \multicolumn{9}{c|}{$\theta$} \\ 
			\hline
			$\delta$ & 0.1 & 0.2 & 0.3 & 0.4 & 0.5 & 0.6 & 0.7 & 0.8 & 0.9 \\ 
			0 & (1.31,18) & (1.38,17) & (1.53,15) & (1.61,14) & (1.82,12) & (2.10,10) & (2.49, 8) & (2.49, 8) & (2.49, 8) \\ 
			0.1 & (1.20,20) & (1.25,19) & (1.37,17) & (1.51,15) & (1.69,13) & (1.93,11) & (2.46, 8) & (2.46, 8) & (2.46, 8) \\ 
			0.2 & (1.10,22) & (1.19,20) & (1.30,18) & (1.43,16) & (1.59,14) & (1.79,12) & (2.22, 9) & (2.42, 8) & (2.42, 8) \\ 
			0.3 & (1.02,24) & (1.10,22) & (1.19,20) & (1.29,18) & (1.49,15) & (1.67,13) & (2.03,10) & (2.39, 8) & (2.39, 8) \\ 
			0.4 & (0.98,25) & (1.05,23) & (1.13,21) & (1.23,19) & (1.41,16) & (1.65,13) & (2.01,10) & (2.36, 8) & (2.36, 8) \\ 
			0.5 & (0.92,27) & (0.98,25) & (1.09,22) & (1.18,20) & (1.34,17) & (1.55,14) & (1.86,11) & (2.34, 8) & (2.34, 8) \\ 
			0.6 & (0.89,28) & (0.95,26) & (1.05,23) & (1.13,21) & (1.27,18) & (1.46,15) & (1.84,11) & (2.31, 8) & (2.31, 8) \\ 
			0.7 & (0.84,30) & (0.92,27) & (1.01,24) & (1.08,22) & (1.21,19) & (1.45,15) & (1.71,12) & (2.29, 8) & (2.29, 8) \\ 
			0.8 & (0.82,31) & (0.89,28) & (0.95,26) & (1.04,23) & (1.21,19) & (1.38,16) & (1.70,12) & (2.27, 8) & (2.27, 8) \\ 
			0.9 & (0.80,32) & (0.87,29) & (0.92,27) & (1.04,23) & (1.16,20) & (1.31,17) & (1.60,13) & (2.25, 8) & (2.25, 8) \\ 
			1 & (0.78,33) & (0.82,31) & (0.89,28) & (1.00,24) & (1.11,21) & (1.30,17) & (1.59,13) & (2.06, 9) & (2.23, 8) \\ 
			\hline
			\multicolumn{10}{|c|}{$T_{double}=40$} \\ 
			\hline
			$T_{lag}=7$ & \multicolumn{9}{c|}{$\theta$} \\ 
			\hline
			$\delta$ & 0.1 & 0.2 & 0.3 & 0.4 & 0.5 & 0.6 & 0.7 & 0.8 & 0.9 \\ 
			0 & (0.85,30) & (0.93,27) & (1.00,25) & (1.07,23) & (1.20,20) & (1.38,17) & (1.61,14) & (2.10,10) & (2.49, 8) \\ 
			0.1 & (0.78,33) & (0.85,30) & (0.90,28) & (0.99,25) & (1.11,22) & (1.25,19) & (1.51,15) & (1.93,11) & (2.46, 8) \\ 
			0.2 & (0.74,35) & (0.78,33) & (0.85,30) & (0.93,27) & (1.02,24) & (1.19,20) & (1.43,16) & (1.79,12) & (2.42, 8) \\ 
			0.3 & (0.74,35) & (0.74,35) & (0.78,33) & (0.87,29) & (0.96,26) & (1.10,22) & (1.29,18) & (1.67,13) & (2.39, 8) \\ 
			0.4 & (0.74,35) & (0.74,35) & (0.74,35) & (0.82,31) & (0.92,27) & (1.05,23) & (1.23,19) & (1.65,13) & (2.36, 8) \\ 
			0.5 & (0.74,35) & (0.74,35) & (0.74,35) & (0.78,33) & (0.87,29) & (0.98,25) & (1.18,20) & (1.55,14) & (2.34, 8) \\ 
			0.6 & (0.74,35) & (0.74,35) & (0.74,35) & (0.74,35) & (0.82,31) & (0.95,26) & (1.13,21) & (1.46,15) & (2.31, 8) \\ 
			0.7 & (0.74,35) & (0.74,35) & (0.74,35) & (0.74,35) & (0.80,32) & (0.92,27) & (1.08,22) & (1.45,15) & (2.29, 8) \\ 
			0.8 & (0.74,35) & (0.74,35) & (0.74,35) & (0.74,35) & (0.76,34) & (0.89,28) & (1.04,23) & (1.38,16) & (2.27, 8) \\ 
			0.9 & (0.74,35) & (0.74,35) & (0.74,35) & (0.74,35) & (0.74,35) & (0.87,29) & (1.04,23) & (1.31,17) & (2.25, 8) \\ 
			1 & (0.74,35) & (0.74,35) & (0.74,35) & (0.74,35) & (0.74,35) & (0.82,31) & (1.00,24) & (1.30,17) & (2.06, 9) \\ 
			\hline
			\multicolumn{10}{|c|}{$T_{double}=50$} \\ 
			\hline
			$T_{lag}=7$ & \multicolumn{9}{c|}{$\theta$} \\ 
			\hline
			$\delta$ & 0.1 & 0.2 & 0.3 & 0.4 & 0.5 & 0.6 & 0.7 & 0.8 & 0.9 \\ 
			0 & (0.76,34) & (0.80,32) & (0.88,29) & (0.96,26) & (1.07,23) & (1.20,20) & (1.45,16) & (1.82,12) & (2.49, 8) \\ 
			0.1 & (0.74,35) & (0.74,35) & (0.80,32) & (0.87,29) & (0.96,26) & (1.11,22) & (1.31,18) & (1.69,13) & (2.46, 8) \\ 
			0.2 & (0.74,35) & (0.74,35) & (0.74,35) & (0.80,32) & (0.90,28) & (1.02,24) & (1.24,19) & (1.59,14) & (2.42, 8) \\ 
			0.3 & (0.74,35) & (0.74,35) & (0.74,35) & (0.76,34) & (0.85,30) & (0.96,26) & (1.14,21) & (1.49,15) & (2.39, 8) \\ 
			0.4 & (0.74,35) & (0.74,35) & (0.74,35) & (0.74,35) & (0.80,32) & (0.92,27) & (1.09,22) & (1.41,16) & (2.17, 9) \\ 
			0.5 & (0.74,35) & (0.74,35) & (0.74,35) & (0.74,35) & (0.76,34) & (0.87,29) & (1.05,23) & (1.34,17) & (2.15, 9) \\ 
			0.6 & (0.74,35) & (0.74,35) & (0.74,35) & (0.74,35) & (0.74,35) & (0.82,31) & (0.98,25) & (1.27,18) & (1.97,10) \\ 
			0.7 & (0.74,35) & (0.74,35) & (0.74,35) & (0.74,35) & (0.74,35) & (0.80,32) & (0.95,26) & (1.21,19) & (1.96,10) \\ 
			0.8 & (0.74,35) & (0.74,35) & (0.74,35) & (0.74,35) & (0.74,35) & (0.76,34) & (0.92,27) & (1.21,19) & (1.94,10) \\ 
			0.9 & (0.74,35) & (0.74,35) & (0.74,35) & (0.74,35) & (0.74,35) & (0.74,35) & (0.89,28) & (1.16,20) & (1.80,11) \\ 
			1 & (0.74,35) & (0.74,35) & (0.74,35) & (0.74,35) & (0.74,35) & (0.74,35) & (0.86,29) & (1.11,21) & (1.79,11) \\ 
			\hline
			\multicolumn{10}{|c|}{$T_{double}=80$} \\ 
			\hline
			$T_{lag}=7$ & \multicolumn{9}{c|}{$\theta$} \\ 
			\hline
			$\delta$ & 0.1 & 0.2 & 0.3 & 0.4 & 0.5 & 0.6 & 0.7 & 0.8 & 0.9 \\ 
			0 & (0.58,46) & (0.62,43) & (0.67,39) & (0.72,36) & (0.80,32) & (0.93,27) & (1.07,23) & (1.38,17) & (2.10,10) \\ 
			0.1 & (0.66,40) & (0.69,38) & (0.74,35) & (0.74,35) & (0.74,35) & (0.85,30) & (0.99,25) & (1.25,19) & (1.93,11) \\ 
			0.2 & (0.74,35) & (0.74,35) & (0.74,35) & (0.74,35) & (0.74,35) & (0.78,33) & (0.93,27) & (1.19,20) & (1.79,12) \\ 
			0.3 & (0.74,35) & (0.74,35) & (0.74,35) & (0.74,35) & (0.74,35) & (0.74,35) & (0.87,29) & (1.10,22) & (1.67,13) \\ 
			0.4 & (0.74,35) & (0.74,35) & (0.74,35) & (0.74,35) & (0.74,35) & (0.74,35) & (0.82,31) & (1.05,23) & (1.65,13) \\ 
			0.5 & (0.74,35) & (0.74,35) & (0.74,35) & (0.74,35) & (0.74,35) & (0.74,35) & (0.78,33) & (0.98,25) & (1.55,14) \\ 
			0.6 & (0.74,35) & (0.74,35) & (0.74,35) & (0.74,35) & (0.74,35) & (0.74,35) & (0.74,35) & (0.95,26) & (1.46,15) \\ 
			0.7 & (0.74,35) & (0.74,35) & (0.74,35) & (0.74,35) & (0.74,35) & (0.74,35) & (0.74,35) & (0.92,27) & (1.45,15) \\ 
			0.8 & (0.74,35) & (0.74,35) & (0.74,35) & (0.74,35) & (0.74,35) & (0.74,35) & (0.74,35) & (0.89,28) & (1.38,16) \\ 
			0.9 & (0.74,35) & (0.74,35) & (0.74,35) & (0.74,35) & (0.74,35) & (0.74,35) & (0.74,35) & (0.87,29) & (1.31,17) \\ 
			1 & (0.74,35) & (0.74,35) & (0.74,35) & (0.74,35) & (0.74,35) & (0.74,35) & (0.74,35) & (0.82,31) & (1.30,17) \\ 
			\hline
	
	\caption{Optimal dosage and optimal number of fractions ($d^*,N^*$) in the presence of uncertainty in tumor parameters (denoted by $\theta$) and uncertainty in $\rho$ parameters (denoted by $\delta$).}
	\label{table:tumoruncertainty}
\end{longtable}
\end{landscape}

%\begin{figure}[!h]
%	\centering
%	\includegraphics[width=5in]{Histogram4}
%	\caption{Histogram of the price of robustness for all 400 simulations.}
%	\label{fig:histogram}
%\end{figure}

%\begin{figure}[!h]
%	\centering
%	\includegraphics[width=5in]{HistInsFrac}
%	\caption{Infeasibility of the nominal solution when $\rho$ varies inside the uncertainty interval.}
%	\label{fig:histInsFrac}
%\end{figure}

%\begin{figure}[!h]
%	\centering
%	\includegraphics[width=6in]{HistOutFrac}
%	\caption{Infeasibility when $\rho$ varies outside the uncertainty interval: the left panel is for nominal, the right panel is for robust.}
%	\label{fig:histOutFrac}
%\end{figure}

\end{document}